\documentclass[12pt]{amsart}
\usepackage{comment,amsfonts,amssymb,amsmath,amscd}
\let\germ=\mathfrak
\let\gg=\Gamma

\let\cal=\mathcal
\newcommand{\Res}{\mathrm{Res }}
\newcommand{\Ind}{\mathrm{Ind }}
\newcommand{\hgt}{\mathrm{ht}}
\newcommand{\rank}{\mathrm{rank}}
\newtheorem{theorem}{\bf Theorem}[subsection]
\newtheorem{corollary}[theorem]{\bf Corollary}
\newtheorem{lemma}[theorem]{\bf Lemma}
\newtheorem{definition}[theorem]{\bf Definition}
\newtheorem{proposition}[theorem]{\bf Proposition}
\numberwithin{equation}{section}
\title[A notion of rank for unitary representations]
{A notion of rank for unitary representations of reductive groups
based on Kirillov's orbit method}
\author{Hadi Salmasian}
\date{August 28 2005}
\address{Queen's University\\
Department of Mathematics and Statistics\\
Jeffery Hall, University Avenue\\
Kingston, Ontario K7L 3N6\\
Canada
}
\email{hadi@mast.queensu.ca}
\begin{document}
\maketitle
\begin{abstract}
We introduce a new notion of rank for unitary
representations of semisimple groups over a
local field of characteristic zero. The theory is
based on Kirillov's method of orbits for nilpotent
groups over local fields. When the semisimple group is
a classical group, we prove that the new theory
is essentially equivalent to Howe's theory of
$N$-rank \cite{Howe}, \cite{LiInv}, \cite{Scara}.
Therefore, our results provide a systematic
generalization of the notion of a small representation
(in the sense of Howe)
to exceptional groups.
However, unlike previous works which used 
ad-hoc methods to study different types of 
classical groups (and some exceptional ones \cite{Marty},
\cite{savinloke}),
our definition is simultaneously applicable to
both classical and exceptional groups. 
The most important result of 
this paper is a general ``purity'' result
for unitary representations, which demonstrates
how similar partial results in the previous authors' works 
should be formulated  
and proved for an arbitrary
semisimple group
in the language of Kirillov's theory.
 The purity result is a crucial 
step towards
studying small 
representations of exceptional groups.
New results concerning 
small unitary representations of exceptional 
groups will be published in a forthcoming paper \cite{hadi2}.\\
\end{abstract}
\setcounter{tocdepth}{1}
\tableofcontents

\section{Introduction}

\subsection{Background}

The theory of $N$-rank
\cite{Howe}, \cite{LiInv}, \cite{Scara} 
provides an 
effective tool for studying 
``small" unitary representations
of classical semisimple groups. It is based on 
analysing the
restriction of a unitary representation to large
commutative subgroups of such groups.
The study of the notion of $N$-rank
began in \cite{Howe} where the theory
was developed for
the metaplectic group. Roughly
speaking, one can attach equivalence
classes of bilinear symmetric forms to unitary
representations, and these equivalence
classes
canonically correspond to orbits of the adjoint
action of the Levi component of the Siegel
parabolic on its unipotent radical. Later on, the theory
was extended by J. S. Li \cite{LiInv}
and R.Scaramuzzi \cite{Scara} to classical
groups, and an intrinsic connection with the
theta correspondence was also found in
\cite{LiComp}. Representations
which correspond to degenerate classes are
called {\it singular} representations.
They form a large class of non-tempered
representations and therefore their study
is essentially complementary to that
of the tempered representations.
It turns out that this class of representations
is quite well-behaved. In fact Howe 
gave a construction of a large class of such 
representations
\cite{howeconstruct}, and they are applied in
construction
of automorphic forms
\cite{Howeautomorphic},\cite{Liautomorphic}.\\

The lack of a theory of $N$-rank which applies 
to the exceptional groups is in a sense a defect of
this theory due to substantial interest in the
small representations of these groups. 
Aside from some generalizations in the works
of
Loke and Savin \cite{savinloke}
and Sahi, most recently
Weissman \cite{Marty} provides an analogue
for $p$-adic split simply-laced
groups excluding $E_8$ using
the representation-theoretic analogue of
a Fourier-Jacobi functor. \\

Here we develop
a theory of rank which can also be applied
to (almost) all exceptional
groups. The word ``almost'' will be 
explained more precisely in the remark after
Proposition
\ref{complexheisenbergparabolicstructure}.
Our
theory relies on the structure of unipotent
radicals of
certain parabolic subgroups which are naturally
expressible
as a sequence of extensions by Heisenberg
groups. In this paper
these unipotent radicals
are called
{\it H-tower} groups
(see Definition \ref{OKPdefinition}).\\

The main point of this work is to show that although
the abelian nilradicals in
terms of which rank was defined for classical
groups are usually not available in
exceptional cases, the maximal unipotent subgroups
of all semisimple groups have a common feature which
enables us to define a similar notion of rank, which
is essentially equivalent to the older one
for classical groups. This result, apart from
being interesting itself,
provides the main 
tool needed to
study small representations of exceptional groups.
One interesting application of this result is 
that it enables us to obtain sharp bounds on
the behaviour of matrix coefficients of 
unitary representations of exceptional 
groups, a problem studied in \cite{LiZhu}, \cite{Oh}
and \cite{savinloke}.
Results along this line of research will be 
published
elsewhere \cite{hadi2}.\\

It is also possible to generalize Theorem 4.2
of \cite{Howe}, which is the next step towards 
classification of small representations of 
exceptional groups. However, there does not seem
to be any strong connection between the ``exceptional
theta correspondences'' which result from the 
author's work and the ones that already exist in the 
literature (see \cite{SAHDEV}, \cite{MAG} or \cite{Ginz}
for instance). This makes the classification problem 
both harder and more interesting at the same time.
This is a work still in its preliminary stages.\\

This manuscript is organized as follows. Section 
\ref{notation}
is devoted to the notation used throughout the article.
In section 2 we recall
some basics about the Heisenberg groups and the 
Weil representation.
In section 3 we define a specific 
unipotent subgroup of a simple
group with regard to which our rank is defined. 
The ``rankable'' representations of this unipotent
group are defined and studied briefly in section 4. 
Section 5 is devoted to the statement and 
proof of the main theorems, which essentially 
assert that rank is 
a nontrivial invariant
of a representation.
Section 6 shows why our new theory 
generalizes the older one.

\subsection*{Acknowledgement} This paper is part of
the author's thesis written under Roger Howe's 
supervision.
The author thanks him for his support and
encouragement. The author thanks
Igor Frenkel, Gopal Prasad, Gordan Savin,
Siddhartha Sahi,
Jeb Willenbring and Gregg Zuckerman for fruitful
conversations/communications, and especially 
the anonymous referee for providing
the author with several extremely
helpful suggestions and for mentioning the 
inaccuracies which existed in section 5
of an earlier version of the paper.

\subsection{Notation}\label{notation}
In this section we introduce some
notation which will be used
throughout this article. Throughout this
paper, we will be working with a local field
$\mathbb F$ which is either $\mathbb R,\mathbb C$ or
a finite extension of $\mathbb Q_p$, $p\neq 2$.
Let $\overline{\mathbb F}$ be the algebraic closure of $\mathbb F$,
and
$\mathbf{G}$ be the $\overline{\mathbb F}$-points of an absolutely simple,
simply connected linear algebraic group
defined over $\mathbb F$, and
$\mathbf{G}_{\mathbb F}$ be the
group of $\mathbb F$-rational points of
$\mathbf G$. We assume $\mathbf G$ is 
$\mathbb F$-isotropic.
By the Kneser-Tits conjecture, which holds in
the case of local fields
(see Proposition 7.6 and Theorem 7.6 in
\cite[\S 7.2]{plra}), 
$\mathbf G_\mathbb F$ 
is equal to 
the group generated by its unipotent elements.
(In the archimedean case, 
$\mathbf G_\mathbb F$ is a connected Lie group.)
Let $G$ be a finite topological
central extension of 
$\mathbf G_\mathbb F$.
(For the definition of a topological central
extension, the reader is referred to
\cite[Chapter 1]{mooreextension}.)\\

Take a maximally split Cartan subgroup
$\mathbf H$ of $\mathbf G$ (which is also
defined
over $\mathbb F$), and let
$\mathbf A$ be the maximal split
torus inside $\mathbf H$.
Let $A$ be the
inverse image of $\mathbf A_{\mathbb F}$ in
$G$.\\

We denote the Lie algebras of the groups $\mathbf G,
\mathbf H, \mathbf A, G, A$ by
$\germ{g}, \germ{h},
\germ{a}$, $\germ{g}_\mathbb F,
\germ{a}_\mathbb F$
respectively.
More precisely, $\germ g, \germ h, \germ a$ are
$\overline{\mathbb F}$-Lie algebras and 
$\germ g_\mathbb F$
and $\germ a_\mathbb F$ are $\mathbb F$-forms of the
corresponding algebras.\\

Let $\Delta$ be an absolute root 
system associated to
$\germ h$, and  let $\Sigma$ be
the restricted (or relative, as some people call it)
root system associated to $\germ a$.
The root space corresponding to any
$\alpha\in\Delta$ is denoted by $\germ g_\alpha$,
and the coroot for $\alpha$ 
is denoted by $H_\alpha$.\\

We fix a $\mathbf G$-invariant symmetric bilinear form 
$(\cdot,\cdot)$ on $\germ g$, 
normalized in such a way that
if for any two roots $\beta,\gamma$ 
we define $(\beta,\gamma)$ to be $(H_\beta,H_\gamma)$, then
$(\beta,\beta)=2$
for any long root $\beta$. All roots of
simply-laced root systems are considered as long ones.
For
$
\germ h_\mathbb F
=\germ h\cap \germ g_\mathbb F
$
we have
$$
\germ h_\mathbb F=\germ a_\mathbb F\oplus
\germ t_\mathbb F
$$
where $\germ t_\mathbb F$
is the orthogonal complement of 
$\germ a_\mathbb F$ in $\germ h_\mathbb F$
with respect to the 
invariant symmetric bilinear form 
$(\cdot,\cdot)$ introduced above.\\

As usual, $\Delta^+$ and $\Sigma^+$ denote
the set of positive roots. We assume 
$\Delta^+$, when restricted to
$\germ a$, contains $\Sigma^+$.
Let
$\Delta_B$ and $\Sigma_B$ denote the bases
for positive roots in the
corresponding root systems such that 
${\Delta_B}$, when restricted
to ${\germ a}$, contains $\Sigma_B$.\\

The choice of a positive system determines
a fixed Borel subalgebra $\germ b$ of $\germ g$
(or a minimal parabolic ${\germ p}_{\min}$
of $\germ g_\mathbb F$).
A parabolic subalgebra of $\germ g$
(respectively, of $\germ g_\mathbb F$)
is called {\it standard}
if it contains $\germ b$ (respectively, 
$\germ p_{\min}$).
We identify the standard parabolic 
subgroups and subalgebras of $\germ g$ 
(respectively of $\germ g_\mathbb F$)
by subsets
of  $\Delta_B$ (respectively of 
$\Sigma_B$)
they naturally correspond to.
More precisely,
for any subset $\Phi\subseteq  \Delta_B$,
by $\germ p_\Phi$ we
mean the standard parabolic subalgebra
of $\germ g$ which contains all root spaces
$\germ g_\alpha$ for any $\alpha\in\Delta^+$
as well as any $\germ g_{-\alpha'}$
such that
$\alpha'\in\Delta^+$ is in the semigroup
generated by $\Phi$.
The definition will be similar for
the parabolic subalgebra (this time
of $\germ g_\mathbb F$)
corresponding to any $\Phi\subseteq\Sigma_B$.
Later on, when there is no risk of confusion
between the parabolics of $\germ g$ and $\germ g_\mathbb F$,
we use the more concise notation
$\germ p_\Phi$ for both situations.
Also, again for simplicity 
and when there is no risk of ambiguity,
Lie subalgebras of both $\germ g$ and 
$\germ g_\mathbb F$ may be
denoted by
Gothic letters in a similar way; i.e. 
without a subscript 
``$\mathbb F$'' in the case of 
a Lie subalgebra of $\germ g_\mathbb F$.\\

The standard parabolic subgroups of
$\mathbf G$ and $G$ are denoted by
$\mathbf P_\Phi$ and $P_\Phi$ respectively.
Note that a parabolic subgroup of $G$
means the inverse image of a parabolic
subgroup of $\mathbf G_\mathbb F$ in $G$.
In other words, if the short exact sequence
\begin{equation}
\label{centralextensionofG}
\begin{CD}
1@>>> F@>i>>G@>p>>\mathbf G_\mathbb F@>>> 1
\end{CD}
\end{equation}
(where $F$ is a finite subgroup of the center of $G$)
represents our topological central extension,
then a parabolic $P$ of $G$ is of the form
$P=p^{-1}(\mathbf P_\mathbb F)$ for an
$\mathbb F$-parabolic $\mathbf P$ of 
$\mathbf G$.
It follows from 
\cite[II 11 Lemme]{dufloacta} (or
from an explicit construction of the 
universal central extension by
Deodhar \cite[Section 1.9]{deod}) that  
topological central extensions of semisimple
groups over local fields are 
split over the unipotent subgroups
\footnote{In fact an abstract (and hence topological) 
universal
central extension of $\mathbf G_\mathbb F$ splits 
over the maximal
unipotent subgroup.
I am indebted to Gopal Prasad for 
pointing to reference \cite{deod}.}
and therefore
one can have an analogous  
``Levi'' decomposition in $G$ as well;
i.e. if $P=p^{-1}(\mathbf P_\mathbb F)$
and $\mathbf P=\mathbf L\mathbf N$
is a Levi decomposition of 
$\mathbf P$ (with $\mathbf L$
being the reductive factor and $\mathbf N$
being the unipotent radical)
then  
one can express $P$ as
$P=LN$,
where:

\begin{enumerate}

\item[1.]
$L=p^{-1}(\mathbf L_\mathbb F)$.\\

\item[2.] $N\subseteq p^{-1}(\mathbf N_\mathbb F)$
is a normal subgroup of $P$
and the map 
$$p:N\mapsto \mathbf N_\mathbb F$$
is an isomorphism of topological groups.
(Existence of $N$ is what follows from the
results of Duflo or Deodhar cited above.)

\end{enumerate}

We may as well drop the subscript $\Phi$ of
the standard parabolics 
$P_\Phi$ and $\mathbf P_\Phi$
when there is
no risk of confusion about the identifying
subset $\Phi$ of simple positive roots.
Moreover, whenever a parabolic 
subgroup or subalgebra is standard,
we assume that the Levi decomposition considered
above is also standard.\\

Let $K$ be an arbitrary locally compact group.
Then the center of
$K$ will be denoted by $\mathcal Z(K)$.
Unitary representations of $K$ 
are defined in the usual way (see \cite{Mackey})
and denoted by Greek letters.
The unitary dual of $K$ (which is the set of 
equivalence classes of the irreducible unitary
representations of $K$ endowed with the Fell 
topology) is denoted by $\hat{K}$. Now suppose 
$K'$ is a closed subgroup of $K$ and let 
$\sigma,\sigma'$ be unitary representations of $K,K'$
respectively. Then, as usual, $\mathrm{Res}_{K'}^K\sigma$ and
$\mathrm{Ind}_{K'}^K\sigma'$ mean restriction and 
unitary induction. If there is no ambiguity about $K$, we 
use the simpler notation $\sigma_{|K'}$ instead of
$\mathrm{Res}_{K'}^K\sigma$. \\

Suppose  $\germ l$ is an arbitrary Lie algebra.
Then
the center
of $\germ l$ is denoted by $\germ z(\germ l)$.\\


\section{Heisenberg group
and Weil representation}
\label{theheisenberggroupandtheweilrepresentation}

\subsection{Heisenberg groups and algebras}
\label{theheisenberggroup}

Let $n$ be a positive integer. 
An $\mathbb F$-Heisenberg group
(or simply a Heisenberg group)
$H_n$ is a two-step
nilpotent
$2n+1$-dimensional $\mathbb F-$group
which is isomorphic to the subgroup
of $GL_{n+2}(\mathbb F)$ of unipotent, upper-triangular 
matrices which do not have 
nonzero elements outside the diagonal,
the first row and the last column.
This group has
a one-dimensional
center. We denote the Lie algebra
of $H_n$ by $\germ h_n$.
Since $H_n\subset GL_{n+2}(\mathbb F)$,
$\germ h_n$ sits in the usual way inside 
$\germ{gl}_{n+2}(\mathbb F)$.\\

The Lie algebra $\germ h_n$ of 
$H_n$ can also 
be described more abstractly 
as 
$$\germ h_n=W_n\oplus \germ z$$
with the Lie bracket introduced below,
where $W_n$ is a $2n$-dimensional 
space endowed
with a nondegenerate symplectic form 
$<\cdot,\cdot>$, and $\germ z$ is the center 
of $\germ h_n$, a one-dimensional Lie subalgebra.
We fix a Lie algebra isomorphism between 
$\germ z$ and 
$\mathbb F$, and denote the element in $\germ z$
corresponding to
$1\in \mathbb F$  by $Z$.
The Lie bracket on $W_n$ is defined as follows: 
\begin{eqnarray*}
    \textrm{for every }X,Y\in W_n\ ,\qquad [X,Y]=
    <\!\!X,Y\!\!>\!\!Z.
\end{eqnarray*}

Viewing $\germ h_n$ as a subalgebra of 
$\germ{gl}_{n+2}(\mathbb F)$,
the exponential map of $GL_{n+2}(\mathbb F)$
gives a bijection between 
$\germ h_n$ and $H_n$, and 
$\mathcal Z(H_n)$ is equal to
$\exp(\germ z)$. This bijection is the 
exponential map of $H_n$.

\subsection{Stone-von Neumann theorem}
\label{heisenbergrepresentations}

We briefly mention 
the structure of irreducible, 
infinite-dimensional unitary
representations of a Heisenberg 
group. For a more 
detailed discussion, the reader can study
\cite{HoweBulletin}, \cite{CG} or \cite[Chapter 1]{taylor} 
for instance.\\

The Stone-von Neumann theorem states
that for any nontrivial
unitary character of $\mathcal Z(H_n)$, 
up to unitary equivalence there is a unique 
infinite-dimensional irreducible unitary
representation of 
$H_n$ having this central 
character. We will describe an explicit realization 
(the so-called Schr\"{o}dinger
model) for this representation below.\\

Consider an arbitrary polarization
of $W_n$, i.e. a decomposition of
$W_n$ as a direct sum
\begin{equation}
    \label{polarizationofw}
W_n=\germ x_n\oplus\germ y_n
\end{equation}
of maximal isotropic subspaces.
It is possible to choose 
bases $\{X_1,\ldots, X_n\}$ and
$\{Y_1,\ldots, Y_n\}$ of $\germ x_n$ and
$\germ y_n$ respectively
such that 
\begin{equation}
    \label{normality}
<X_i,Y_j>=\delta_{i,j}.
\end{equation}
There is an
isomorphim 
\begin{equation}
\label{dualanditself}
e:\germ x_n\approx\germ y_n^*
\end{equation}
obtained via the symplectic form.
$\mathcal Z(H_n)\approx\exp(\mathbb F)$ and any
nontrivial multiplicative (unitary) 
character of $\mathcal Z(H_n)$
is equal to
$\chi\circ\exp^{-1}$, for some $\chi$,
where $\chi$ can be  
any nontrivial 
additive (unitary) 
character of $\germ z(\germ h_n)\approx \mathbb F$.
Such characters $\chi$ are in
one-to-one correspondence with 
elements of $\mathbb F-\{0\}$.
For example,
if $\mathbb F=\mathbb R$, 
$\chi$ will be of the form 
$\chi(x)=e^{\mathrm ixy}$ 
for some $y\in\mathbb R-\{0\}$ (here 
$\mathrm i=\sqrt{-1}$).
This follows from the well-known fact that 
the unitary
dual of a local field is identifiable to 
itself (see \cite{WEILNUM}).\\

The (irreducible) representation $\rho$ of $H_n$ 
with central character $\chi$
may be realized on the Hilbert space
$$
{\mathcal H}_{\rho}=L^2(\germ y_n)
$$
as follows. Let $y=\sum_{i=1}^ny_iY_i\in \germ y_n$.
Then $H_n$ acts on
${\mathcal H_{\rho}}$ as  
\begin{eqnarray}
   (\rho(e^{tX_i})f)(y)&=&\chi
(-ty_i)f(y)\notag\\
    (\rho(e^{tY_i})f)(y)&=&f(y+tY_i)\\
    \label{heisenbergaction}
    (\rho(e^{tZ})f)(y)&=&\chi(t)f(y).\notag
\end{eqnarray}    
This representation is irreducible and 
unitary with respect to
the usual inner product of 
$\mathcal H_{\rho}$.\\

\noindent{\bf Notation.}\ Henceforth, 
whenever a locally compact
group $N$ is isomorphic to
a Heisenberg group, we denote the subset of 
$\hat N$ consisting of 
the family of the infinite-dimenional
representations constructed above 
by ${\hat{N}}_\circ$.

\subsection{The Weil representation} 
\label{thefockmodel}
Let $Sp(W_n)$ be the symplectic
group associated to $W_n$. When
$\mathbb F=\mathbb C$, take $Mp(W_n)=Sp(W_n)$,
and otherwise let
$Mp(W_n)$ be the (metaplectic) 
double covering of $Sp(W_n)$
(see \cite{weil}). 
$Mp(W_n)$ acts through 
$Sp(W_n)$ on 
$H_n\approx\exp(\germ h_n)$
(or equivalently on $\germ h_n$)
as follows. 
For any $w\oplus z\in W_n\oplus\germ z$,
$$
g:w\oplus z\rightarrow g\cdot w\oplus z.
$$
One can consider the 
semidirect product $Mp(W_n)\ltimes
H_n$. We have (see \cite{weil}):

\begin{proposition}
    \label{weilrepresentationextension}
    Let $\chi$ be any nontrivial additive 
character of \ $\mathbb F$.
Then
the irreducible unitary representation $\rho$ of
    $H_n$ with central character $\chi$
can be extended to a unitary 
    representation (still denoted by $\rho$)
    of 
    $Mp(W_n)\ltimes H_n$;
    i.e. for any $g\in Mp(W_n)$ and any $h\in H_n$:
    $$
       \rho(g)\rho(h)
       \rho(g^{-1})=
    \rho(ghg^{-1}).
    $$
\end{proposition}


\section{Construction of H-tower groups}

\subsection{The Heisenberg parabolic}
\label{subsecheisenbergparabolic}
We use the notation of section \ref{notation}.
Let $\tilde{\beta}$ be the highest root in $\Delta$
with respect to the positive system
chosen in section \ref{notation} and 
let $H_{\tilde\beta}$
be the coroot for $\tilde\beta$.
Let
$(\cdot,\cdot)$ denote the 
$\mathbf G$-invariant
symmetric bilinear form introduced in
that section; i.e. 
$(\tilde\beta,\tilde\beta)=2$. 
One obtains a 
grading \[\germ g=
\bigoplus_{j\in\mathbb Z}\germ g_j\]
where $\germ g_j$ is the $j$-eigenspace 
of $\mathrm{ad}({H_{\tilde\beta}})$; i.e.
$$
\germ g_j=\{X\in\germ g : 
[H_{\tilde\beta},X]=jX\}.
$$
Note that in fact
$\germ g_k=\{0\}$
for $|k|>2.$
The Jacobi identity implies  
$[\germ g_i,\germ g_j]\subseteq\germ g_{i+j}$.
Note that
$\germ g_1=\germ g_{-1}=\{0\}$
if and only if $\germ g=\germ{sl}_2$.\\

In the next proposition and 
at some points of this  
paper,
we will consider  
$\germ g$ as above such that:
\begin{equation}
\label{gdefine}
\diamond\ \ \ \
\germ g\neq\germ{sl}_2\textrm{ and }
\tilde\beta\textrm{ is defined over }\mathbb F.
\end{equation} 
Condition (\ref{gdefine}) is not
a crucial assumption. In fact the author believes that
it can be removed. See the discussion in the remark
after Proposition
\ref{complexheisenbergparabolicstructure}.

\begin{proposition}
    \label{complexheisenbergparabolicstructure}
    Let $\germ g$ satisfy condition (\ref{gdefine}). Then
    
    \begin{itemize}
    \item[1.]
    $\germ g_0\oplus\germ g_1\oplus\germ g_2$ 
is a  
parabolic subalgebra,
$\germ  g_1\oplus \germ g_2$ is its nilradical,
and 
$$
[[\germ g_0,\germ g_0],\germ g_2]=\{0\}.
$$
The nilradical is a Heisenberg Lie 
algebra with center $\germ g_2$ and symplectic 
space $\germ g_1$. One can describe the symplectic 
form $<\!\!\cdot,\cdot\!\!>_1$ 
on $\germ g_1$ using
the Lie bracket 
as follows: fix $X_{2}\in
\germ g_{2}-\{0\}$ and
let $<\!\!\cdot,\cdot\!\!>_1$ be 
such that 
\begin{displaymath}
\textrm{for every } X,Y\in \germ g_1\ , 
\quad [X,Y]=<\!\!X,Y
\!\!>_1\!
\cdot X_{2}\ .
\end{displaymath}
Then $<\!\!\cdot,\cdot\!\!>_1$ 
will be nondegenerate, and the 
(adjoint)
representation of $[\germ g_0,\germ g_0]$
on $\germ g_1$ will be 
a symplectic representation
with respect to $<\!\!\cdot,\cdot\!\!>_1$.\\

\item[2.]
When $\germ g$ is not of type $\mathbf A_l$
($l>1$), 
the parabolic of part 1 
is characterized
by the subset 
$\Delta_B-\{\beta\}$ of $\Delta_B$ 
for the unique simple root 
$\beta\in\Delta_B$
which satisfies $(\beta,\tilde{\beta})\neq 0$. 
For type $\mathbf A_l$ ($l>1$), there are two such
simple roots $\beta',\beta''$, 
and the parabolic corresponds to
$\Delta_B-\{\beta',\beta''\}$.\\

\item[3.]
When $\germ g$ is not of type $\mathbf A_l$ ($l>1$), 
$\beta$ 
corresponds to a simple restricted root (still denoted
by $\beta$).
For $\germ g$ of type $\mathbf A_l$ ($l>2$)
and $\germ g_\mathbb F$ non-split,
the pair $\{\beta',\beta''\}$ corresponds to a 
simple restricted 
root $\beta$ 
in the restricted root system. \\
\end{itemize}

\noindent Let the sets $S\subset \Delta_B$ and $T\subset \Sigma_B$
be defined as follows.
Set $S=\{\beta\}$ when $\germ g$ is not
of type $\mathbf A_l$ ($l>1$) and 
$S=\{\beta',\beta''\}$ otherwise.
Also,
set $T=\{\beta',\beta''\}$ when $\germ g_\mathbb F$
is split and $\germ g$ is of type 
$\mathbf A_l$ ($l>1$), and 
$T=\{\beta\}$ otherwise. \\

\begin{itemize}

\item[4.]
The parabolic 
$\germ p_{\Delta_B-S}$ is defined over $\mathbb F$,
and $(\germ p_{\Delta_B-S})_\mathbb F
=\germ p_{\Sigma_B-T}$.
Its nilradical
is a Heisenberg algebra.
$\germ p_{\Delta_B-S}$ 
is a maximal parabolic 
when $\germ g$ is not of type
$\mathbf A_l$, ($l>1$).
\end{itemize}

\end{proposition}

\begin{proof}
    See \cite[\S 2]{GW} and 
\cite[\S 10]{torassoduke}.
\end{proof}
{\noindent \bf Remark.} 
The assumption that $\germ g$ 
satisfies condition (\ref{gdefine})
holds for all but a very small class
of Lie algebras $\germ g$.
Fox example, when
$\mathbb F=\mathbb R$, 
the cases where
$\tilde\beta$ is not defined over $\mathbb R$
are $\germ f_{4(-20)}$, $\germ e_{6(-26)}$, 
$\germ{sp}(p,q)$, and $\germ{sl}(n,\mathbb H)$. 
The notion of rank in
this paper is based on the existence of 
a unipotent
subgroup of $G$ which is expressible
as a tower of extensions by
Heisenberg groups (see Definition
\ref{OKPdefinition}), and therefore it is not
applicable to the groups mentioned above.
However, in all of 
these groups there does exist
a similar structure which is 
called an OKP subgroup 
in 
\cite{HoweOKP}. Therefore in principle 
the
main results of this paper (especially 
Theorem \ref{main})
should generalize
to
groups associated with these real forms.
However, addressing the technical 
problems which
arise with including those cases in this
paper makes it much more technical.
To keep our presentation as simple and 
uniform as possible, we
do not include those special cases in this paper.
The author intends to deal with those cases
elsewhere.

\begin{definition}
    \label{heisenbergparabolicdefinition}
    The 
    parabolic subalgebras 
    $\germ p_{\Delta_B-S}$ of
    $\germ g$ and 
    $\germ p_{\Sigma_B-T}$ of $\germ g_\mathbb F$
    or their corresponding parabolic
    subgroups 
    $\mathbf P_{\Delta_B-S}$ of $\mathbf G$
    and
    $P_{\Sigma_B-T}=
    (\mathbf P_{\Delta_B-S})_\mathbb F$ 
    of $\mathbf G_\mathbb F$ (or the parabolic of $G$
which corresponds to $P_{\Sigma_B-T}$)
are called the 
    {\it Heisenberg}
    parabolics. 
\end{definition}

    We tend to drop their identifying 
    subscripts for simplicity when there is no risk 
    of confusion.\\

Table I below demonstrates the structure of
$\germ g_1$ as a representation of 
$[\germ g_0,\germ g_0]$.
Here $V_\varpi$ denotes the representation with
highest weight $\varpi$ and  
$\varpi_i$ denotes the $i$-th 
fundamental weight of the corresponding 
Lie algebra. 
We use the notation of 
\cite[Planche I]{bourbaki} 
for numbering fundamental weights.
See also \cite{Tits} for more explicit
information. \\

\begin{center}
\begin{tabular}{|l|c|c|}
    \hline
    $\germ g$&$[\germ g_0,\germ g_0]$&$\germ g_1$\\
    \hline
    $\mathbf{A_l}\  (\mathbf l\geq 3)$&
    $\mathbf{A_{l-2}}$&$V_{\varpi_1}\oplus
    V_{\varpi_1}^*$\\
    $\mathbf{B_2}$&$\mathbf{A_1}$&
    $V_{\varpi_1}$\\
    $\mathbf{B_3}$&$\mathbf{A_1}\times\mathbf{A_1}$&
    $V_{\varpi_1}\hat\otimes V_{\varpi_1}$\\
    $\mathbf{B_l}\  (\mathbf l\geq 4)$&
    $\mathbf{A_1}\times\mathbf{B_{l-2}}$&
    $V_{\varpi_1}\hat\otimes V_{\varpi_1}$\\
$\mathbf{C_2}$&$\mathbf{A_1}$&$V_{\varpi_1}$\\    
$\mathbf{C_l}\ (\mathbf l\geq 3)$&
    $\mathbf{C_{l-1}}$&$V_{\varpi_1}$\\
    $\mathbf{D_4}$&$\mathbf{A_1}\times
    \mathbf{A_1}\times
	\mathbf{A_1}$&
    $V_{\varpi_1}\hat\otimes V_{\varpi_1}
	\hat\otimes
	V_{\varpi_1}$\\
    $\mathbf{D_5}$&$\mathbf{A_1}\times\mathbf{A_3}$&
    $V_{\varpi_1}\hat\otimes V_{\varpi_2}$\\
    $\mathbf{D_l}\ (\mathbf l\geq 6)$&
	$\mathbf{A_1}\times
    \mathbf{D_{l-2}}$&$V_{\varpi_1}
	\hat\otimes V_{\varpi_1}$\\
    $\mathbf{E_6}$&$\mathbf{A_5}$&$V_{\varpi_3}$\\
    $\mathbf{E_7}$&$\mathbf{D_6}$&$V_{\varpi_6}$\\
    $\mathbf{E_8}$&$\mathbf{E_7}$&$V_{\varpi_7}$\\
    $\mathbf{F_4}$&$\mathbf{C_3}$&$V_{\varpi_3}$\\
    $\mathbf{G_2}$&$\mathbf{A_1}$&$V_{3\varpi_1}$\\
    \hline
    \end{tabular}
   
   \vspace{3mm}
   {\bf Table I}\\
    
\end{center}    
    
\subsection{The H-tower subgroup $N_\gg$ of $G$}
\label{sokpnilradical}
Let $\germ g$ be as in condition 
(\ref{gdefine}).
In this section we describe a 
nilpotent subalgebra of $\germ g$ (and
another one of $\germ g_\mathbb F$) which 
is
fundamental to our definition of rank.
The construction of these nilpotent subalgebras 
is based on what is usually referred to as
Kostant's Cascade.
We show that this nilpotent Lie subalgebra
is actually the nilradical of a parabolic 
subalgebra (see Propositions \ref{ngammatower}
and \ref{structureofcomplexifiedrealngamma}).
The unipotent subgroup of $G$ which
corresponds to this nilpotent subalgebra 
will play an important role in the rest 
of the paper.\\

First assume that $\germ g$ 
splits over $\mathbb F$; 
i.e. all roots of $\germ g$ 
are defined over $\mathbb F$.
Let $\germ p=\germ l\oplus\germ n$ 
be the Heisenberg parabolic of $\germ g$,
obtained by Proposition
\ref{complexheisenbergparabolicstructure},
with the usual Levi decomposition; i.e.
$\germ l$ is the Levi factor and $\germ n$ is
the nilradical of $\germ p$.
In the case of orthogonal algebras  
the commutator $[\germ l,\germ l]$
is not a simple Lie algebra.
In fact when $\germ g$ is
of types $\mathbf B_l\ (l\geq 3)$
or $\mathbf D_l\ (l>4)$, the
commutator $[\germ l,\germ l]$
is a direct sum of
the form
\begin{equation}
\label{llsl2}
[\germ l,\germ l]
=\germ{sl}_2\oplus\germ s
\end{equation}
where $\germ s$ is simple.
When $\germ g$ is of type $\mathbf D_4$,
we have 
$$
[\germ l,\germ l]=
\germ{sl}_2\oplus\germ{sl}_2\oplus
\germ{sl}_2.
$$ 

\begin{definition}\label{MDEFINEM}

Let $\germ m$ be defined as follows:
\begin{itemize}
\item[I.] If $\germ g$ is of types 
$\mathbf D_l\ (l>4)$
or $\mathbf B_l\ (l>3)$ then 
$\germ m$ is equal to the summand $\germ s$
given in
(\ref{llsl2}).\\

\item[II.] If $\germ g$ is of types $\mathbf A_2,\mathbf A_3,
\mathbf B_2=\mathbf C_2,
\mathbf B_3,
\mathbf D_4$ or $\mathbf G_2$ then
$\germ m=\{0\}$.\\

\item[III.] Otherwise, 
$\germ m=[\germ l,\germ l]$.

\end{itemize}

\end{definition}

One can repeatedly apply Proposition 
\ref{complexheisenbergparabolicstructure}
as follows. 
First we apply it to $\germ m$.
If $\germ m$ is nonzero, then 
$\germ m\neq\germ{sl}_2$ and
Proposition
\ref{complexheisenbergparabolicstructure}
guarantees the existence of a Heisenberg
parabolic in $\germ m$. Let $\germ p'$ be this 
parabolic of $\germ m$, 
and let $\germ l'$ be the Levi factor 
of $\germ p'$. Let $\germ m'$ be
the subalgebra of $[\germ l',\germ l']$
which is defined in the same way that $\germ m$ was
defined as a subalgebra of $[\germ l,\germ l]$
in Definition \ref{MDEFINEM}.
Then we apply Proposition 
\ref{complexheisenbergparabolicstructure}
to 
$\germ m'$, 
and so on. This process can be repeated as long as 
Proposition
\ref{complexheisenbergparabolicstructure}
can be applied.
As a result, we
obtain a sequence 
$S_1,\ldots,S_r$ of subsets
of $\Delta_B$, 
where each $S_i$, 
defined as in part 2 of Proposition
\ref{complexheisenbergparabolicstructure}, 
contains either a simple root or a pair of
simple roots.
Each $S_i$ corresponds to 
the ``highest root'' 
$\tilde{\beta}_i$ 
obtained in the $i$-th step.
We also denote the
sequence of 
Heisenberg parabolics by 
$\germ p^1,\ldots,\germ p^r$,
with the Levi decomposition
\begin{equation}
\label{sokpsokpsokp}
\germ p^j=\germ l^j\oplus\germ n^j.
\end{equation}
Therefore, $\germ p^1$ is the Heisenberg
parabolic of $\germ g$, $\germ p^2$ is 
the Heisenberg parabolic of $\germ m$,
and so on.
Each $\germ n^j$ is normalized by any
$\germ l^{j'}$ for $j'\geq j$ and hence by
$\germ n^{j'}$. For a similar reason,
$\germ n^{j'}$ acts trivially on the center of
$\germ n^j$.
Each $\germ n^j$ is isomorphic to
a Heisenberg algebra $\germ h_{d_i}$
for some $d_i$.
Therefore the Lie algebra
$\germ n^1\oplus\cdots\oplus\germ n^r$ 
is a tower of successive 
extensions by 
Heisenberg $\overline{\mathbb F}$-algebras. 
The following
proposition is obvious.

\begin{proposition}
    \label{complexngammadefinition}
    $\germ n^1\oplus\cdots\oplus\germ n^r$ is 
    equal to the nilradical of the 
    parabolic 
    subalgebra $\germ p_\mathbf \gg$ of $\germ g$
    where
    
    $$
    \mathbf \gg=\Delta_B-(S_1\cup\ldots\cup S_r).
    $$
\end{proposition}

Now assume $\germ g$ satisfies
condition 
(\ref{gdefine}) but it is not necessarily
$\mathbb F$-split.
Recall that the positive system
for $\Sigma$ is compatible with the one 
for $\Delta$. 
We can apply
Proposition
\ref{complexheisenbergparabolicstructure}
repeatedly again.
Note that the number of possible iterations 
for a non-split $\mathbb F$-form of 
$\germ g$ is often smaller than the number of
possible iterations in the split case.
(In certain groups with small rank, they may 
be equal.) This is because the successive
Levi factors may fail to satisfy condition 
(\ref{gdefine}).\\

Part 3 of Proposition
\ref{complexheisenbergparabolicstructure}
yields a sequence
$T_1,\ldots,T_s$ of subsets of
$\Sigma_B$ (for some $s\leq r$)
where each $T_j$ 
contains a simple restricted root or a pair of them.
Again we obtain a nested sequence
$$
\germ p^1\supset\cdots\supset\germ p^s
$$ 
of Heisenberg parabolics (this time inside
$\germ g_\mathbb F$) 
and the nilradical of any $\germ p^i$ 
is normalized by the nilradical of 
any $\germ p^{i'}$ when $i<i'$.
Therefore $\germ p^1$ is the Heisenberg
parabolic subalgebra of $\germ g_\mathbb F$. 
We warn the reader that the new $\germ p^i$'s
are different from those which appear
immediately before Proposition \ref{complexngammadefinition};
in fact the older $\germ p^i$'s are obtained from the newer
$\germ p^i$'s by an extension of scalars.\\

\newpage

\addtolength{\evensidemargin}{-1.5cm}
\addtolength{\textwidth}{2.5cm}

\begin{center}
\begin{tabular}{|lr|l|c||lr|l|c|}
\hline
$\germ g_\mathbb F$&&$\germ m_\mathbb F$&$s$&
$\germ g_\mathbb F$&&$\germ m_\mathbb F$&$s$ \\
\hline
$A_{2,2}$&				& -- &1&
$A_{3,3}$&				&--&1\\
$A_{r,r}$&$(r\geq 4)$ 			&$A_{r-2,r-2}$&$\lfloor{r\over 2}\rfloor$&
$\ ^2A_{3,2}^{(1)}$& 			&--&1\\
$\ ^2A_{2r-1,r}^{(1)}$&$(r\geq 3)$	& $\ ^2A_{2r-3,r-1}^{(1)}$&$r-1$&
$\ ^2A_{2,1}^{(1)}$& 			& --&1 \\
$\ ^2A_{2r,r}^{(1)}$&$(r\geq 2)$  	& $\ ^2A_{2r-2,r-1}^{(1)}$&$r$&
$\ ^2A_{3,1}^{(1)}$& 			& --&1\\
$\ ^2A_{2r+1,r}^{(1)}$&$(r\geq 2)$ 	& $\ ^2A_{2r-1,r-1}^{(1)}$&$r$&
$B_{3,3}$& 				& --&1\\
$B_{4,4}$& 				& $C_{2,2}$&2&
$B_{r,r}$&$(r\geq 4)$ 			& $B_{r-2,r-2}$&$\lfloor{ r\over 2}\rfloor$\\
$B_{3,2}$& 				& --&1&
$B_{4,3}$& 				& --&1\\
$B_{r,r-1}$&$(r\geq 5)$ 		& $B_{r-2,r-3}$&$\lfloor {r-1\over 2}\rfloor$&
$C_{2,2}$& 				& --&1\\
$C_{r,r}$&$(r\geq 3)$ 			& $C_{r-1,r-1}$&$r-1$&
$\ ^1D_{4,4}^{(1)}$& 			& --&1\\
$\ ^1D_{5,5}^{(1)}$&			& $A_{3,3}$&2&
$\ ^1D_{r,r}^{(1)}$&$(r\geq 6)$ 	& $\ ^1D_{r-2,r-2}^{(1)}$&$\lfloor{r-1\over 2}\rfloor$\\
$\ ^1D_{4,2}^{(1)}$&			& --&1&
$\ ^1D_{5,3}^{(1)}$& 			& --&1\\
$\ ^1D_{r+2,r}^{(1)}$&$(r\geq 4)$ 	& $\ ^1D_{r,r-2}^{(1)}$&$\lfloor {r\over 2}\rfloor$&
$\ ^2D_{4,3}^{(1)}$& 			& --&1\\
$\ ^2D_{5,4}^{(1)}$&			&$\ ^2 A_{3,2}^{(1)}$&2&
$\ ^2D_{r+1,r}^{(1)}$&$(r\geq 5)$ 	& $\ ^2D_{r-1,r-2}^{(1)}$&$\lfloor 
{r-1\over 2}\rfloor$\\
$\ ^1D_{4,2}^{(2)}$ &				& -- &1&
$\ ^1D_{2r,r}^{(2)}$&$(r\geq 3)$ 	& $\ ^1D_{2r-2,r-1}^{(2)}$&$r-1$\\
$\ ^1D_{5,1}^{(2)}$&			& --&1&
$\ ^1D_{2r+3,r}^{(2)}$&$(r\geq 2)$	& $\ ^1D_{2r+1,r-1}^{(2)}$&$r$\\
$\ ^2D_{5,2}^{(2)}$&			& $\ ^2A_{3,1}^{(1)}$&2&
$\ ^2D_{2r+1,r}^{(2)}$&$(r\geq 3)$	& $\ ^2D_{2r-1,r-1}^{(2)}$&$r$\\
$\ ^2D_{4,1}^{(2)}$&			&--&1&
$\ ^2D_{2r+2,r}^{(2)}$&$(r\geq 3)$	& $\ ^2D_{2r,r-1}^{(2)}$&$r$\\
$\ ^3D_{4,2}^2$	&			& --&1&
$\ ^6D_{4,2}^2$	&			& --&1\\
$\ ^1E_{6,2}^{16}$&			& --&1&
$\ ^1E_{6,6}^0$	&			& $A_{5,5}$&3\\
$\ ^2E_{6,4}^2$	&			& $\ ^2A_{5,3}^{(1)}$&3&
$E_{7,4}^9$	&			& $\ ^1D_{6,3}^{(2)}$&3\\
$E_{7,7}^0$	&			& $\ ^1D_{6,6}^{(1)}$&3&
$E_{8,8}^0$	&			& $E_{7,7}^0$&4 \\
$F_{4,4}^0$	&			& $C_{3,3}$&3&
$G_{2,2}^0$	&			& --&1\\
\hline
\end{tabular}

\vspace*{.3cm}
{\small\bf Table II}

\end{center}

\begin{center}
\footnotesize
\begin{tabular}{|lr|l|c||lr|l|c|}

\hline
$\germ g_\mathbb F$ && $\germ m_\mathbb F$&$s$&
$\germ g_\mathbb F$ && $\germ m_\mathbb F$&$s$\\

\hline

$\germ{sl}_3(\mathbb R)$&       	&--				&1&
$\germ{sl}_4(\mathbb R)$&		&--				&1\\
$\germ{sl}_n(\mathbb R)$&$(n\geq 5)$&$\germ{sl}_{n-2}(\mathbb R)$	&$\lfloor {n-1\over 2}\rfloor$&
$\germ{su}(1,q)$	&$(q>1)$	&--				&1\\
$\germ{su}(2,2)$	&		&--				&1&
$\germ{su}(r,q)$	&$(2\leq r<q)$	&$\germ{su}(r-1,q-1)$		&$r$\\
$\germ{su}(q,q)$	&$(3\leq q)$	&$\germ{su}(q-1,q-1)$		&$q-1$&
$\germ{so}(1,q)$	&$(q\geq 4)$	&--				&1\\
$\germ{so}(2,q)$	&$(q\geq 3)$	&--				&1&
$\germ{so}(3,q)$	&$(q\geq 3)$	&--				&1\\
$\germ{so}(4,4)$	&		&--				&1&
$\germ{so}(4,q)$	&$(q\geq 5)$	&   $\germ{so}(2,q-2)$		&2\\
$\germ{so}(r,q)$	&$(5\leq r<q)$  & $\germ{so}(r-2,q-2)$	        &$\lfloor {r\over 2}\rfloor$&
$\germ{so}(q,q)$	&$(5\leq q)$    & $\germ{so}(q-2,q-2)$	        &$\lfloor {q-1\over 2}\rfloor$\\
$\germ{sp}_4(\mathbb R)$&		&--				&1&
$\germ{sp}_{2n}(\mathbb R)$&$(n\geq 3)$	& $\germ{sp}_{2n-2}(\mathbb R)$ &$n-1$\\
$\germ{so}^*(6)$	&		&--				&1&
$\germ{so}^*(8)$	&		&--				&1\\
$\germ{so}^*(2r)$	&$(r\geq 5)$	&$\germ{so}^*(2r-4)$		&$\lfloor {r-1\over 2}\rfloor$&
$(\germ e_6,\germ{sp}_4)$&		&$\germ{sl}_6(\mathbb R)$	&3\\
$(\germ e_6,\germ{su}_6\times \germ{su}_2)$&&$\germ{su}(3,3)$		&3&
$(\germ e_6,\germ{so}(10)\times \germ u(1))$&&$\germ{su}(1,5)$		&2\\
$(\germ e_7,\germ{su}_8)$&		&$\germ{so}(6,6)$		&3&
$(\germ e_7,\germ{so}(12)\times \germ{su}_2)$&&$\germ{so}^*(12)$        &3\\
$(\germ e_7,\germ e_6\times\germ u(1))$&& $\germ{so}(2,10)$		&2&
$(\germ e_8,\germ{so}(12))$&&$(\germ e_7,\germ{su}_8)$			&4\\
$(\germ e_8,\germ e_7\times \germ {su}_2)$&&$(\germ e_7,\germ e_6
\times \germ u(1))$&3&
$(\germ f_4,\germ{sp}_3\times\germ{su}_2)$&&
$\germ{sp}_6(\mathbb R)$						&3\\
$(\germ g_2,\germ{su}_2\times\germ{su}_2)$&&--				&1&&&&\\

\hline

\end{tabular}

\vspace*{3mm}
{\small \bf Table III}
\end{center}

\newpage
\addtolength{\evensidemargin}{1.5cm}
\addtolength{\textwidth}{-2.5cm}

Tables II and III explain how
Proposition \ref{complexheisenbergparabolicstructure}
is applied iteratively to $\germ g_\mathbb F$. 
For simplicity, the real and the $p$-adic cases are 
separated: Table II is for the $p$-adic case and Table III
is for the real case.\\

The 
column $\germ g_\mathbb F$ in Table II shows the Tits index of 
$\germ g_\mathbb F$ and
the column $\germ m_\mathbb F$ 
shows the Tits index of $\germ m_\mathbb F$.
For any $\germ g_\mathbb F$ the number $s$ is given too. 
Real exceptional Lie algebras in Table III
are identified by the symmetric pairs 
$(\germ g_\mathbb R,\germ k_\mathbb R)$.
The only $\mathbb F$-forms of $\germ g$ which appear in the  
columns $\germ g_\mathbb F$ of both
of the tables
are those for which $\germ g$ satisfies condition
(\ref{gdefine}).
If an entry in the column corresponding to $\germ m_\mathbb F$
is equal to ``--'', this means either that the Lie algebra
$\germ m$ is equal to $\{0\}$ or that the Lie algebra
$\germ m$ is semisimple but it does not satisfy condition 
(\ref{gdefine}), i.e.
the highest root of $\germ m$ is not defined
over $\mathbb F$.
Therefore one can use Tables II and III
to see how many times Proposition \ref{complexheisenbergparabolicstructure}
can be applied to a particular $\mathbb F$-form of $\germ g$. \\

We denote the nilradical of the 
parabolic $\germ p^j$ 
of $\germ g_\mathbb F$
by $\germ h^j$.
It is an $\mathbb F$-Heisenberg algebra.\\

\begin{proposition}
    \label{ngammatower}
    The nilpotent Lie subalgebra 
    $\germ h^1\oplus\cdots\oplus\germ h^s$
    of $\germ g_\mathbb F$ 
    is equal to the nilradical of the
    parabolic subalgebra 
    $\germ p_\gg$ of $\germ g_\mathbb F$
    where 
    $$
    \gg=\Sigma_B-(T_1\cup\ldots\cup T_s).
    $$
   
\end{proposition}

Moreover, we have the following proposition, which
describes 
the relationship between $\germ n^i$'s and $\germ h^i$'s.

\begin{proposition}
\label{structureofcomplexifiedrealngamma}    
Let $\germ p_\gg$ be the
parabolic subalgebra 
of $\germ g_\mathbb F$ 
defined in Proposition
\ref{ngammatower}
with Levi decomposition
$\germ p_\gg=\germ l_\gg\oplus\germ n_\gg$.
Then
$$
\germ n_\gg\otimes\overline{\mathbb F}=\germ n^1\oplus
\cdots\oplus\germ n^s
$$
(where the $\germ n^j$'s are the same as those
which appear in Proposition \ref{complexngammadefinition})
and thus 
$\germ n_\gg\otimes{\overline{\mathbb F}}$
is equal to the nilradical of the
parabolic $\germ p_{\mathbf \gg'}$ of $\germ g$
where
$$
{\mathbf \gg'}=\Delta_B-(S_1\cup\ldots\cup S_s).
$$

\end{proposition}

Therefore $\mathbf P_{\mathbf\gg'}$ is an
$\mathbb F$-parabolic of $\mathbf G$
and $(\mathbf P_{\mathbf\gg'})_\mathbb F=P_\gg$, 
where 
$P_\gg$ is the parabolic
of $G$ 
which is associated to 
the set $\gg\subseteq\Sigma_B$
defined in Proposition \ref{ngammatower}.      
$P_\gg$ will play a significant role in
the rest of the paper. Its
Levi decomposition can be written as 
\begin{equation}
    \label{levidecompositionofpgamma}
    P_\gg=L_\gg N_\gg.
\end{equation}

\begin{definition}\label{OKPdefinition}
Let $U$ be the group of $\mathbb F$-points of 
a unipotent linear algebraic
$\mathbb F$-group $\mathbf U$.
$U$ is said to be an ``H-tower'' group 
if and only if  
it satisfies one of the properties {\rm I} or 	
{\rm II} below.
\begin{itemize}
\item[I.] ${U}=\{1\}$; i.e. $U$ is trivial.
\item[II.] $U$ is isomorphic to 
a semidirect product $U'\ltimes U''$ where 
\begin{itemize}
\item[i.]$U''$ is 
an $\mathbb F$-Heisenberg group and is the group of 
$\mathbb F$-points of 
an algebraic $\mathbb F$-subgroup $\mathbf U''$
of 
$\mathbf U$.
\item[ii.] $U'$ is the 
group of $\mathbb F$-rational points
of 
an algebraic subgroup $\mathbf U'$ of $\mathbf U$
and  
the action of $U'$ on  
$U''$ in the semidirect product 
$U'\ltimes U''$
comes from an 
algebraic action (defined over $\mathbb F$)
of 
$\mathbf U'$ on $\mathbf U''$
by group
automorphisms of $\mathbf U''$ 
(which,
a fortiori, leave elements of   
$\mathcal Z(\mathbf U'')$
invariant and act on the symplectic space
$\mathbf U''/\mathcal{Z}(\mathbf U'')$ 
via symplectic operators)
\footnote{This means that there is an 
$\mathbb F$-homomorphism of algebraic groups
$\Phi:\mathbf U'\mapsto\mathbf{Sp}\left(
\mathbf U''/
\mathcal Z(\mathbf U'')\right)$ and consequently 
$\Phi(U')\subseteq Sp\left(U''/\mathcal Z(U'')\right)$ 
.}.
\item[iii.] $U'$ is an H-tower group.
\end{itemize}
\end{itemize}

\end{definition}

Therefore an H-tower group $U$ can be expressed
as a tower of successive extensions:

\begin{equation}\label{OKPDecomp}
U=H^1\cdot
H^2\cdots H^t=H^1\rtimes(\cdots\rtimes(
H^{t-1}\rtimes H^t)\cdots)
\end{equation}
where each $H^j$ is a Heisenberg group; i.e. 
$H^j=H_{d_j}$ for some $d_j$ (see section 
\ref{theheisenberggroup}).\\ 

{\noindent \bf Remarks.} 

\noindent 1. The name 
``H-tower'' is chosen so that
it reminds the reader that 
the group is a tower of extensions by Heisenberg
groups.\\

\noindent 2. Real H-tower groups form a subclass of
OKP groups defined in \cite{HoweOKP}.\\

\noindent 3. The number $t$ in 
(\ref{OKPDecomp}) is referred to as the 
{\it height} of the tower of extensions. 
It will be denoted by $\hgt(U)$. \\

\noindent 4. Consider an H-tower group $U$ which is 
expressed as in (\ref{OKPDecomp}). 
for any $j\in\{1,...,\hgt(U)\}$ we 
denote the quotient group
$H^j\cdots H^{\hgt(U)}\approx U/H^1
\cdots H^{j-1}$ by $U_j$. $U_j$ is also an H-tower group.\\

If in (\ref{OKPDecomp}) 
we take  $U=N_\gg$ and  
$H^j=\exp \germ h^j$, then
by Proposition
\ref{ngammatower} we have
\begin{proposition}
$N_\gg$ is an H-tower group and 
$\hgt(N_\gg)=s$, where $s$ is as in
Proposition
\ref{ngammatower}.

\end{proposition}

\section{Rankable representations of H-tower
groups}

\subsection{Oscillator extension and 
rankable representations}

Fix an H-tower group $U$, expressed as a
tower of successive 
extensions as in (\ref{OKPDecomp}).
Let $\chi_1$ be a nontrivial additive
character of $\mathbb F$, and  
consider the  
unitary representation $\rho_{1}$ of
$H^1$ with central character $\chi_1$
(see section 
\ref{heisenbergrepresentations}). 
When restricted to the 
inverse image (in the metaplectic group)
of the maximal unipotent 
subgroup of the symplectic group, 
the Weil representation factors through
a representation of the maximal unipotent 
subgroup of the symplectic group.
(In fact the two-fold central extension of 
the symplectic group splits over the maximal 
unipotent subgroup and, more importantly, 
this splitting is unique. This holds because
any two splittings would differ by a 
finite-order character, whereas the maximal 
unipotent subgroup is a divisible group and therefore 
it does not have such characters.)
Therefore, any representation
$\rho$ of any Heisenberg group $H_n$
with central
character $\chi$ 
(see section \ref{heisenbergrepresentations})
is extendable to the unipotent
radical of any Borel subgroup of the 
symplectic group $Sp(W_n)$.
This 
implies that the representation 
$\rho_{1}$ of $H^1$ is
extendable (in at least one way) to $U$. 
We still denote this extension by
$\rho_{1}$. 
Note that the extension of
$\rho_{1}$ is not 
necessarily unique, as
one can for example twist it by a 
one dimensional representation
of $U$ which
is trivial on $H^1$. However, since
the extension of $\rho$ 
to the metaplectic group 
is unique, we can distinguish 
the extension of $\rho_1$ 
obtained by
the restriction of the Weil representation
without ambiguity. Let $\rho_1$ mean this specific
extension.\\

In a similar fashion, for any $j>1$ one can 
construct a representation of 
$$U_j=U/H^1\cdots H^{j-1}\approx H^j\cdots 
H^{\hgt(U)}$$
as follows. 
We take an arbitrary nontrivial
additive character $\chi_j$
of 
$\mathbb F$, and the representation
$\rho_{j}$ of $H^j$
with central character $\chi_j$.
As before, we use the Weil representation
to extend $\rho_j$ 
to $ U_j$. 
Since 
$$
U_j\approx U/ H^1\cdots H^{j-1},
$$ 
this representation can be extended to
$U$, trivially on $H^1\cdots H^{j-1}$. 
Thus using $\rho_j$ 
we have constructed a representation of
$U$.
We still keep
the notation $\rho_{j}$ for this 
representation. 

\begin{definition}\label{s-OKPdef}
Let $U$ be an H-tower group described as
in (\ref{OKPDecomp}). 
A representation of $U$
is called ``rankable'' if and only if it is
unitarily equivalent to a tensor product of the 
form 
$$\rho_{1}\otimes\rho_{2}
\otimes\cdots\otimes\rho_{k}
$$
for some $k\leq\hgt(U)$, 
where $\rho_j$
is the representation of 
$U$ which is obtained (using
the Weil representation) by extending 
the irreducible representation of
$H^j$
with central character $\chi_j$. 
The rank of a $k$-fold tensor product 
(including the case $k=0$,
i.e. the trivial 
representation) is defined to be $k$.
If $\sigma$ is a rankable representation 
of rank $k$, we write $\rank(\sigma)=k$
\end{definition}

\noindent{\bf Remark. }It follows from 
Corollary \ref{OKPORBITSDIM}
that $\mathrm{rank}(\sigma)$ is
well-defined; i.e. it is an invariant of the
unitary equivalence class of $\sigma$.

\subsection{Kirillov theory 
for rankable representations}

We denote the coadjoint orbit attached 
to an irreducible unitary 
representation $\sigma$  of a 
nilpotent group
by $\mathcal O^*_\sigma$.
(See 
\cite{Kirillov}, \cite[Chapter 3]{Kirillov2},
\cite{CG}
or \cite[\S2.2]{AlexanderKirillov}
for elaborated treatments of 
Kirillov's orbital theory.)
In his seminal paper \cite{Kirillov},
Kirillov developed his method of 
orbits for
simply connected nilpotent real 
Lie groups. 
However, later in his 1966 ICM lecture 
he explained that essentially the same 
theory can be applied to algebraic
unipotent groups over $p$-adic fields.
See also \cite[\S4]{moore} for more details. \\

\begin{lemma}\label{constructible}
Let $N$
be the group of $\mathbb F$-rational points of 
a unipotent linear 
algebraic $\mathbb F$-group 
$\mathbf N$.
Let $\sigma$
be an irreducible unitary representation of 
$N$ and let 
$\mathcal O_\sigma^*$ be the coadjoint orbit 
attached to it. Then $\mathcal O_\sigma^*$
is an analytic manifold
in the sense of \cite[Chapter III, \S2]{serre}.

\end{lemma}
\begin{proof}
The lemma follows immediately from 
\cite[\S 3.1, Corollary 2]{plra}.
\end{proof}

Lemma \ref{constructible} implies
that one can speak of 
the
dimension of a coadjoint orbit. 
(The transitivity under the group action
implies that the dimension is the same around 
every point.)
For any coadjoint
orbit  $\mathcal O^*$, let 
$\dim \mathcal O^*$ 
denote its dimension.

\begin{lemma}
    \label{N_1N_2}
    Let $\mathbf N$ be a unipotent algebraic group
such that $\mathbf N=\mathbf N'\ltimes\mathbf N''$ 
as algebraic $\mathbb F$-groups 
(and the action of $\mathbf N'$
on $\mathbf N''$ is defined over $\mathbb F$) 
and let $N=N'\ltimes N''$
be the group of $\mathbb F$-rational 
points of 
$\mathbf N$ (where $N',N''$ are
the $\mathbb F$-points of $\mathbf N',\mathbf N''$).
Let the Lie algebras 
of $N',N''$ be
    $\germ n',\germ n''$ respectively.
    The Lie algebra of the 
    semidirect product $N=N'\ltimes N''$ 
    is $\germ n=\germ n'\oplus\germ n''$ as a 
    vector space, and we have a canonical isomorphism
    of dual spaces 
$\germ n^*={\germ n'}^*\oplus{\germ n''}^*$.
    Let $\sigma',\sigma''$ 
be irreducible
    unitary
    representations of
    $N',N''$ respectively,
such that
\begin{itemize}
\item[1)] $\sigma'$ is an 
irreducible unitary representation of
$N'$,
extended trivially on $N''$ to 
a representation 
$\tilde{\sigma}'$ of $N$.
\item[2)]$\sigma''$ 
is an irreducible unitary representation of $N''$
which extends to a unitary representation 
$\tilde{\sigma}''$ of $N$. (In other words, 
$\mathrm{Res}_{N''}^N\tilde{\sigma}''=\sigma''$.)
\end{itemize}

Then $\tilde\sigma'\otimes\tilde\sigma''$
is an irreducible representation of $N$
and

$$
\mathrm{dim}\ 
\mathcal O^*_{\tilde\sigma'\otimes\tilde
\sigma''}=
\mathrm{dim}\ \mathcal O^*_{\sigma'}
    +\mathrm{dim}\ \mathcal O^*_{\sigma''}.$$

\end{lemma}

\begin{proof}

Consider the map $\mathbf j$ defined as the composition
$$N\mapsto N\times 
N\mapsto N/N'' \times N\approx N'\times N$$
where the leftmost 
map is the diagonal embedding and the
middle map (i.e. the map $N\times N\mapsto N/N''\times N$)
is projection onto the first factor (i.e. 
it is given by $(n_1,n_2)\mapsto (n_1N'',n_2)$). 
Observe that

\begin{equation}\label{localtahdid}
\tilde{\sigma}'\otimes\tilde{\sigma}''=
\mathrm{Res }^{N'\times N}_{\mathbf j(N)}
\sigma'\hat{\otimes}{\tilde{\sigma}''}.
\end{equation}

It is easily seen by  
Mackey theory \cite[Theorem 3.12]{Mackey} 
that 
$\tilde{\sigma}'
\otimes\tilde{\sigma}''$ is 
irreducible. Now take a maximal 
chain of analytic subgroups 
$$
N'\times N=N^0\supset N^1\supset N^2\supset\cdots
\supset \mathbf j(N)
$$ 
such that each $N^j$ has codimension 1 in 
$N^{j-1}$. It follows from (\ref{localtahdid})
that if $\sigma'\otimes\tilde{\sigma}''$ 
is considered as a representation of $N^{j-1}$,
then 
$$\mathrm{Res }^{N^{j-1}}_{N^j}
\sigma'\otimes\tilde{\sigma}''$$
is an irreducible representation of $N^j$.
Now it follows from the statement
of Theorem 2.5.1 of
\cite{CG}
that for any such $j$ 
we are in the situation of part (b) 
of that theorem. But then Proposition 1.3.4 of
\cite{CG}
implies that in this case, the dimension of 
the coadjoint orbit of 
$\sigma'\otimes\tilde{\sigma}''$ 
considered as
a representation of $N^j$ is the same as the
dimension of the coadjoint orbit of 
$\sigma'\otimes\tilde{\sigma}''$ 
considered as a representation 
of $N^{j-1}$. Consequently, 
 $$\mathrm{dim}\ \mathcal O^*_{\sigma'\hat{\otimes}
    \tilde{\sigma}''}=
    \mathrm{dim}\  \mathcal O^*_{\tilde{\sigma}'
    \otimes\tilde{\sigma}''}.$$

Now consider the 
map $\mathbf q$ such that 
$$\mathbf q:N'\times N''\mapsto N'\times N$$
which is defined to be the identity
map $N'\mapsto N'$ in the first component and 
the injection $N''\subset N$
in the second component.
Then 
$$
\mathrm{Res}_{\mathbf q(N'\times N'')}^{N'\times N}
\sigma'\hat{\otimes}\tilde{\sigma}''=
\sigma'\hat{\otimes}
\sigma''
$$
which is an irreducible representation.
Again a similar argument (taking a 
maximal chain, etc.)
proves that 
$$\mathrm{dim}\ \mathcal O^*_{
\sigma'\hat\otimes \tilde\sigma''}=\mathrm{dim}
\ \mathcal O^*_{\sigma'
\hat{\otimes}\sigma''}.
$$
But the coadjoint orbit
$\mathcal O^*_{\sigma'
\hat{\otimes}\sigma''}$ 
equals
$\mathcal O^*_{\sigma'}\times
\mathcal O^*_{\sigma''}$,
which completes the proof. 
\end{proof}

One can apply Lemma \ref{N_1N_2} 
to any H-tower group iteratively
and obtain the following result.

\begin{corollary}
    \label{OKPORBITSDIM}
Let $U$ be an H-tower group as in
(\ref{OKPDecomp}). Let 
$$\sigma=\rho_{1}\otimes\rho_{2}
\otimes\cdots\otimes\rho_{k}$$
be a rankable representation 
of $U$.
(Recall that $\rho_j$ is extended
from a 
representation of
$H^j$ with central character $\chi_j$.) 
Then $\sigma$ is irreducible. Moreover, 
as an analytic manifold, the coadjoint
orbit $\mathcal O^*_\sigma$ 
has dimension 
$2(n_1+\cdots+n_k)$,
where 
$$
n_j={\mathrm{dim} \germ h^j-1\over 2}.
$$
\end{corollary}

\section{Main theorems}

\subsection{Technical issues of central extensions}

\label{mackeyanalysissection}

Let $G,\mathbf G, \germ g$ be as in section
\ref{notation} such that $\germ g$
satisfies condition (\ref{gdefine}).
Let $P$ be the Heisenberg parabolic of 
$G$ (see Definition 
\ref{heisenbergparabolicdefinition}).
Suppose $P=LN$ is the Levi decomposition
of $P$. Therefore $N$ is a Heisenberg group. 
As in section
\ref{heisenbergrepresentations}
let $\rho$ be the irreducible unitary representation of 
$N$ with an arbitrary 
nontrivial central character $\chi$.
Let $N_\gg$ be the H-tower subgroup of $G$ 
constructed in section \ref{sokpnilradical}
and assume $\hgt(N_\gg)>1$.\\

This section is mainly devoted to a technical issue which
arises with central extensions of $G$. The main result is
Proposition \ref{definingp1}. 
To avoid being distracted from 
the main point of the paper, the reader may assume
Proposition \ref{definingp1}
and go on to section \ref{makeyII}.\\

Proposition 
\ref{complexheisenbergparabolicstructure}
implies that the adjoint action provides 
a group homomorphism $[L,L]\mapsto Sp(N/\mathcal Z(N))$. 
Our next aim is to somehow apply 
Proposition \ref{weilrepresentationextension} 
and extend 
$\rho$ to a representation of a larger 
group which at least contains $N_\gg$. 
One natural candidate could be 
$[P,P]$.
However, when 
$\mathbb F\neq\mathbb C$,
in order 
to extend the representation $\rho$ 
of $N$ to $[P,P]$ we need to know that indeed  
the group homomorphism  
$[L,L]\mapsto Sp(N/\mathcal Z(N))$ 
is a composition of a group homomorphism
into the {\it metaplectic} group, i.e. 
$$[L,L]\mapsto Mp(N/\mathcal Z(N)),$$ 
and 
the projection map 
$$Mp(N/\mathcal Z(N))\mapsto Sp(N/\mathcal Z(N)).$$
This issue 
was also dealt with in 
\cite{torassoduke}
(see \cite[Section 3.3]{GanSavin}).
We will explain it more clearly below.
In fact it can be seen that
in some cases, if 
$\mathbf P_\mathbb F=\mathbf L_\mathbb F
\mathbf N_\mathbb F$
is the Heisenberg parabolic of 
$\mathbf G_\mathbb F$,
then 
$[\mathbf L_\mathbb F,\mathbf L_\mathbb F]$
does
{\it not} act as a 
subgroup
of the metaplectic group, and therefore
extending $\rho$ to 
$[\mathbf P_\mathbb F,\mathbf P_\mathbb F]$ 
may not 
be possible. An obvious 
example is 
when $\mathbf G_\mathbb F$ is the 
symplectic group. 
A less obvious example is the split group of type 
$\mathbf F_4$. (Actually,
in both cases the extension is possible
if we use the metaplectic covers. This
may be seen after some simple calculations.
For real groups
one can use the results of 
\cite{AdamsTrapa} which tell us when
a covering of a real simple Lie group
splits over embedded root subgroups 
$SL_2(\mathbb R)$. However, here we do
not rely on such calculations.)\\

To overcome this difficulty, we will
see below that we have to 
consider
an appropriate finite (topological)
central extension of $G$ instead of $G$ 
itself. Note that a 
representation of $G$ can trivially be
considered as a representation of its 
central extension, and we can always study
this extension instead of the original $G$
because the central extension
has an H-tower subgroup too, 
which is identical to $N_\gg$  
and is essentially obtained by exponentiating
the corresponding nilpotent subalgebra 
$\germ n_\gg$
of $\germ g_\mathbb F$. This follows 
from what was explained in section \ref{notation} 
as well,
that a universal topological central 
extension of $\mathbf G_\mathbb F$
splits over the maximal unipotent subgroup
(see \cite[II 11 Lemme]{dufloacta} or \cite[Section 1.9]{deod}).
(In fact the splitting is unique because
any two splittings would differ by a 
finite-order character, whereas any 
unipotent subgroup is a divisible group and therefore 
it does not have such characters.) 
Therefore, 
as far as it concerns the main results of
this work,  we can substitute 
$G$ with such a ``good'' central extension of $G$, without
any loss of generality.\\

To show the existence of this ``good'' central extension
for an arbitrary $G$, first we show its existence for
the group $\mathbf G_\mathbb F$. That is to say, we
show that we can find
a 
finite topological central extension 
$\tilde{\mathbf G}_\mathbb F$ of 
$\mathbf G_\mathbb F$ in which 
the representation $\rho$ of $N$
can be extended to a larger subgroup. 
We warn the reader that we may (and will)
think of $N$ as 
a subgroup of both $\mathbf G_\mathbb F$ 
and $\tilde{\mathbf G}_\mathbb F$ 
at the same time, since the latter group
has a subgroup identical to the subgroup $N$ of the former
one.\\

Note that once we find 
$\tilde{\mathbf G}_\mathbb F$,
a simple argument implies that the required
extension exists for $G$ as well.
In fact 
the existence of a universal topological
central extension of $\mathbf G_\mathbb F$
implies that there will be a finite
topological central extension 
which covers both $G$ and 
$\tilde{\mathbf G}_\mathbb F$. This extension
will be the one in which the extension of $\rho$
to a larger subgroup is possible. \\

Next we show that we can find $\tilde{\mathbf G}_\mathbb F$.
Let 
\begin{equation*}
\begin{CD}
1@>>>F@>\hat i>>
\tilde{\mathbf G}_\mathbb F
@>\hat p>>\mathbf G_\mathbb F@>>>1
\end{CD}
\end{equation*}
be a finite topological central extension of 
$\mathbf G_\mathbb F$. 
Let $\mathbf P=\mathbf L\mathbf N$ 
be the Heisenberg
parabolic of $\mathbf G$ with its usual Levi 
decomposition. 
Let $L=\hat p^{-1}({\mathbf L}_{\mathbb F})$
and let $N$ be the unipotent radical of
$P=\hat p^{-1}(\mathbf P_\mathbb F)$.
The group $[\mathbf L,\mathbf L]$, the derived
subgroup of $\mathbf L$, is a product of isotropic and
anisotropic factors (see \cite[\S 8.19]{torassoduke}).
Let $\mathbf L^\mathrm{is}$ be the 
isotropic factor of
$[\mathbf L,\mathbf L]$.
Let $\mathbf M$
be the simply connected subgroup of $[\mathbf L,\mathbf L]$
whose Lie algebra is equal to $\germ m$,
where 
$\germ m$ is the Lie algebra introduced 
in Definition \ref{MDEFINEM}.
(The simply connectedness of 
$\mathbf M$ follows from \cite[\S 8.19]{torassoduke}.)
It follows from 
$\hgt(N_\gg)>1$ that $\mathbf M\subseteq 
\mathbf L^\mathrm{is}$.
Let $M=\hat p^{-1}(\mathbf M_\mathbb F)$.\\

Recall that the adjoint action provides a map 
$$M\mapsto Sp(N/\mathcal Z(N)).$$ 
(See Proposition 
\ref{complexheisenbergparabolicstructure}.) 
Let us call the inverse image of any
subgroup
$M'$ of $M$ inside
the metaplectic group 
$
Mp\left(N/\mathcal Z(N)\right)
$
the {\it metaplectic 
extension} 
of $M'$.
It follows from 
\cite[\S8.25 and \S9.3]{torassoduke} that if
the $\mathbb F$-points of the minimal 
nilpotent $\mathbf G$-orbit of $\germ g$
contain a 
$\tilde{\mathbf G}_\mathbb F$-admissible orbit
(see \cite[\S 3.17]{torassoduke} for the definition
of an admissible orbit), then there
exists a closed subgroup $M'$ of $M$,
which is a normal subgroup of $L$, such that
$\hat p(M')\supseteq \mathbf M_\mathbb F^+$
(where $\mathbf M_\mathbb F^+$ is the subgroup
of $\mathbf M_\mathbb F$ generated by its 
$\mathbb F$-rational unipotent elements)
and
the metaplectic extension 
of $M'$ is trivial. 
(Note that in fact $\hat p(M')\supseteq 
\mathbf M_\mathbb F$, because by 
\cite[\S 8.19]{torassoduke}, $\mathbf M$
is simply connected and therefore 
$\mathbf M_\mathbb F^+=\mathbf M_\mathbb F$
(see \cite[\S 7.2]{plra}).)
This 
means 
that
one can extend
$\rho$ to $M'\cdot N$. 
We keep the notation $\rho$ for such an extension.
We will see below that the extension of 
$\rho$ is actually uniquely
determined on a slightly smaller group, 
and therefore using $\rho$ to denote the 
extension is not 
ambiguous.\\

Tables in the end of \cite{torassoduke} determine
which groups $\mathbf G_\mathbb F$ 
have finite topological
central extensions 
$\tilde{\mathbf G}_\mathbb F$
with
admissible orbits. Among all 
$\mathbb F$-forms
of 
groups 
$\mathbf G$ whose Lie algebras $\germ g$
satisfy condition
(\ref{gdefine}), the only groups which 
do not have coverings with admissible orbits are
groups of Tits index $B_{r,r}$ and $B_{r,r-1}$ 
in the nonarchimedean case and $SU(p,q)$ 
and $SO(p,q)$ when $p+q$ is odd in the 
archimedean case.
In these cases, however, 
the metaplectic extension of 
$\mathbf M_\mathbb F$
is trivial
by \cite[Lemma 2.2]{LiComp} and
\cite[Lemma 7]{kazhdan}. 
Therefore if $M'=\mathbf M_\mathbb F$
then $\rho$ extends to $M'$.
This settles the
issue of existence of an appropriate 
covering for $\mathbf G_\mathbb F$.\\


For a general $G$ we can take $\tilde G$ as a finite 
topological central extension which covers both $G$ and 
$\tilde{\mathbf G}_\mathbb F$. Let $M_0$ be the inverse image of $M'$ in
$\tilde G$.
Obviously extending
$\rho$ to $M_0\cdot N$ is possible.
$M_0$ has a closed finite-index
subgroup $M_1$ which is perfect 
(i.e. $[M_1,M_1]=M_1$).
$M_1$ is found as follows. 
Define the sequence $M^{(0)}=M_0$ and 
$M^{(i)}=[M^{(i-1)},M^{(i-1)}]$. Then each 
$M^{(i)}$ 
is
a subgroup of $M'$ of finite
index and the indices are uniformly bounded
(in fact one upper bound is the size of the kernel 
of the map $\tilde G\mapsto \mathbf G_\mathbb F$). 
Therefore there is some 
$i_\circ$ such that 
$M^{(i_\circ)}=M^{(i_\circ+1)}$; i.e.
$M^{(i\circ)}$ is perfect. Let 
$M_1=\overline{M^{(i_\circ)}}$, the closure
of $M^{(i_\circ)}$ inside $M'$. 
Then $M_1$ is perfect by 
\cite[\S 8.11]{torassoduke}.\\

The 
extension of $\rho$ to $M_1\cdot N$ is 
unique. In fact the projective representation
of $M_1$ is uniquely determined by
the relation
$$\rho(m)\rho(n)\rho(m^{-1})
=\rho(mnm^{-1})\qquad\textrm{for any}
\quad n\in N, m\in M_1
$$
and 
since $M_1$ is perfect, this projective 
representation 
lifts to a linear representation of $M_1$
in a unique way. \\

The discussion in the previous paragraph 
proves the following result.

\begin{proposition}
\label{definingp1}
Let $\mathbf G$ be as in
section \ref{notation} such that 
$\germ g$ satisfies condition
(\ref{gdefine}).
Let $G$ be 
a finite topological central extension of 
$\mathbf G_\mathbb F$ represented as in
(\ref{centralextensionofG}).
Let $P,L,N,\rho$ be as in the beginning of section
\ref{mackeyanalysissection}.
Then there exists a finite topological central 
extension $\tilde G$ of $\mathbf G_\mathbb F$ 
which covers $G$ as in
\begin{equation*}
\begin{CD}
1@>>>\tilde F@>\tilde i>>\tilde G
@>\tilde p>>G@>>>1
\end{CD}
\end{equation*}
and a (closed) subgroup $P_1=M_1\cdot
\tilde N$
of $\tilde G$ such that

\begin{itemize}

\item[i.] $p\circ \tilde p(M_1)=
\mathbf M_\mathbb F$ 
and 
$p\circ \tilde p(\tilde N)
=\mathbf N_\mathbb F$.\\

\item[ii.] $\tilde N$ is the unipotent radical of
$\tilde p^{-1}(P)$; i.e. 
$\tilde p:\tilde N\mapsto N$ is an 
isomorphism.\\

\item[iii.] The representation $\rho$ 
(which can naturally be thought of
as a representation of 
$\tilde N$ too) can be extended in a unique way 
to a representation 
of $P_1$.\\

\item[iv.] $M_1$ is perfect; i.e.
$[M_1,M_1]=M_1$.\\

\item[v.] $P_1$ is a normal subgroup of 
$\tilde p^{-1}(P)$.

\end{itemize}
\end{proposition}

\subsection{Mackey analysis}
\label{makeyII}

In this section we
assume that the notation is the same as
in section \ref{mackeyanalysissection}.
We also
assume (without loss of generality)
that 
\begin{equation}
\label{ASSUMPTION}
\diamond\quad 
G \textrm{\ is equal to
the extension\ } \tilde G
\textrm{\ of  Proposition
\ref{definingp1} and }[G,G]
\textrm{ is dense in }G.
\end{equation}
Note that if $[G,G]$ is dense in $G$, then
from \cite[\S 8.11]{torassoduke}
it follows that in fact $[G,G]=G$.
The assumption $[G,G]=G$ is not a crucial one.
This is because one can always find a 
closed finite-index subgroup $G_1$ of $G$ 
such that $G_1$ is perfect and 
$G=G_1\cdot\mathcal Z(G)$, and once a representation
$\pi$ of $G$ is understood on $G_1$, It is easy to
describe it as a representation of $G$. \\

Assume condition (\ref{ASSUMPTION}) holds. 
Let $P_1\subset G$ be the subgroup 
given in the statement of 
Proposition \ref{definingp1}. 
Then we can write $P_1=M_1\cdot N$, where $M_1$ is
defined in Proposition 
\ref{definingp1} and $N\approx \tilde N$ 
is the unipotent 
radical of the Heisenberg parabolic $P$
of $G$ which appears in parts ii and iii of 
Proposition \ref{definingp1}.\\

Let $\pi$ be a unitary representation of $G$, 
without nonzero
$G$-invariant vectors. 
Consider the restriction of $\pi$
to $P$. 
Recall that $\mathcal Z(N)$ means 
the center of $N$.
By Howe-Moore theorem \cite{HoweMoore},
$\pi$ does not have a nonzero 
$\mathcal Z(N)$-invariant vector 
either. 
Therefore in the direct integral 
decomposition of $\pi$
as a representation of $N$, 
the spectral measure is supported on
$\hat N_\circ$ (see the notation introduced at 
the end
of section \ref{heisenbergrepresentations}).\\

Elements of $P_1$ 
commute
with elements of $\mathcal Z(N)$.
Therefore under the coadjoint action 
of $P$ on the unitary characters of 
$\mathcal Z(N)$, $P_1$ lies
within $\mathrm{stab}_P(\chi)$,
the stabilizer (inside $P$) of any 
nontrivial additive unitary character
$\chi$
of $\mathcal Z(N)$.
In fact all these stabilizers are 
identical groups. 
We denote this common stabilizer group by $J$.\\

The action
of $P$ on $\mathcal Z(N)$ 
(and also on its unitary characters) has 
only a finite number of orbits. 
Consequently, since $\pi$ has no nonzero
$G$-invariant vectors, 
by elementary Mackey theory
(see Theorems 3.11 and 3.12 of \cite{Mackey}) 
the restriction of 
$\pi$ to $P$
can be expressed as a finite direct sum 
$$
\pi_{|P}=\bigoplus_i\mathrm{Ind}_{J}^P\sigma_i
$$
where each $\sigma_i$ is a representation of $J$ 
which, when 
restricted to $N$, is a direct integral of 
representations isomorphic to the representation
$\rho_i$ of $N$ (defined in
section \ref{heisenbergrepresentations})
with some central character $\chi_i$.
However, each $\rho_i$ extends in exactly 
one way
to a 
representation of $P_1$.
 We still denote this extension by $\rho_{i}$.
Again by Mackey theory we can write 
the restriction of $\sigma_i$ to $P_1$
as 
$$
{\sigma_i}_{|P_1}=\nu_i\otimes\rho_{i}
$$
where $\nu_i$ is indeed a representation of $M_1$ 
extended trivially on $N$
to $P_1$.
Therefore we have proved the following result.
  
\begin{lemma}
\label{mackeytheoryforP}
Let $G, P, L$ and $N$ be as in the 
beginning of
section \ref{mackeyanalysissection}
so that $G$ satisfies condition (\ref{ASSUMPTION}).
Let $\pi$
be a unitary representation of $G$ without a 
nonzero $G$-invariant vector.
Then $\pi$ can be
written as a finite direct sum
\begin{equation}
\label{mackeydecompositionofP}
\pi_{|P}=\bigoplus_i 
\mathrm{Ind}_{J}^P\sigma_i
\end{equation}
where each $\sigma_i$ is a representation of $J$
such that 
\begin{equation}
\label{sigmaidecomp}
{\sigma_i}_{|P_1}=\nu_i\otimes \rho_i
\end{equation}
and $\nu_i$ factors to a 
representation
of $M_1$. (In other words, $\nu_i$ 
is trivial on $N$.)
\end{lemma}



Recall the construction of the H-tower subgroup
$U=N_\gg$ (see (\ref{levidecompositionofpgamma})). 
In the notation of (\ref{OKPDecomp}) we have
$H^1=N$. The quotient group 
$U_2=N_\gg/N$ is identical to the 
H-tower subgroup which is constructed for 
$M$. By Proposition \ref{definingp1},
the latter H-tower group should lie within 
$M_1$.

\begin{lemma}\label{rankonedirect}

Let $\rho$ be the extension to $P_1$
of the irreducible representation 
of $N$ (which is also denoted by $\rho$)
with 
an arbitrary 
nontrivial central character $\chi$
(see section \ref{heisenbergrepresentations}).
Let $N_\gg$ be as in 
(\ref{levidecompositionofpgamma}).
Then the restriction of the representation
$\rho$
to $N_\gg$ is supported on rankable 
representations of 
$N_\gg$ of rank one. 
\end{lemma}
\begin{proof}
The uniqueness of 
extension of $\rho$ to
$P_1$ implies that on 
$M_1$, this extension 
should be identical to
the extension obtained by restriction
of extension of $\rho$ 
to the metaplectic 
group acting on 
$N/\mathcal Z(N)$. 
Therefore the lemma follows from
Definition \ref{s-OKPdef}.

\end{proof}


\noindent {\bf Notation.\ }Let $K$ be an arbitrary
(abstract) group and let $K'$ be a subgroup of $K$.
Let $\sigma$ be a representation of $K'$ on a vector
space $\mathcal H$. Let $a\in K$.
Then by $\sigma^a$ we mean a representation of 
the subgroup $K'_a=aK'a^{-1}$ on $\mathcal H$ defined
as follows:
\begin{equation}
\label{Repaction}
\textrm{for every }x\in K'_a\ ,\ 
\sigma^a(x)=\sigma(a^{-1}xa).
\end{equation}
In particular, if $K'=K$, then $\sigma^a$ is
a representation of $K$.\\

At this point we recall Mackey's subgroup theorem.
We will use it in sections 
\ref{proofofthemaintheorems}
and \ref{proofofthemaintheoremsII}.
Giving a precise statement would require
some definitions and would distract us from the 
main point of this paper. Therefore, for a
detailed discussion we refer the reader
to \cite[Theorem 3.5]{Mackey} or 
\cite[Theorem 12.1]{MackeyI}.

\begin{theorem}{\rm (Mackey's Subgroup Theorem)}
\label{MSUBT}
Let $K$ be a locally compact group and let $K',K''$
be its closed subgroups. Assume that $K,K',K''$ 
are ``nice'' (see the above references for details). 
Let $\sigma$ be an irreducible representation
of $K'$. Then the representation
$$\mathrm{Res }^K_{K''} \mathrm{Ind }_{K'}^K\sigma$$
has a direct integral decomposition supported on
representations $\tau_v, v\in K$, where 
$$
\tau_v=\mathrm{Ind }^{K''}_{K''\cap K'_v}
\mathrm{Res }^{K'_v}_{K''\cap K'_v}\sigma^v.
$$ 
Here  
$$
K'_v=vK'v^{-1}
$$ and $\sigma^v$ is defined as in 
(\ref{Repaction}).

\end{theorem}

The following elementary lemma, which essentially
follows from Mackey theory too, 
will help us in 
sections \ref{proofofthemaintheorems}
and \ref{proofofthemaintheoremsII}

\begin{lemma}
\label{orbitinvariance}
Let $K'\subset K$ be two arbitrary 
locally compact
groups such that $K'$ is a closed normal subgroup
of $K$. Let $\sigma$ be a unitary
representation of 
$K$ which acts on the Hilbert space 
$\mathcal H$. Suppose 
$$
\Res^K_{K'}\sigma=\sigma_1\oplus\sigma_2
$$
where $\sigma_1$ and $\sigma_2$ are unitary 
representations of $K'$ such that for any $a\in K$,
we have 
$$
\mathrm{Hom}_{K'}(\sigma_1^a,\sigma_2)=\{0\}.
$$
(Here $\sigma_1^a$ is defined as in (\ref{Repaction})
and 
$\mathrm{Hom}_{K'}(\sigma_1^a,\sigma_2)$
means the space of $K'$-intertwining 
operators
from $\sigma_1^a$ to $\sigma_2$.) 
For $i\in\{1,2\}$,
let $\mathcal H_i$ be the closed 
$\sigma(K')$-invariant subspace of 
$\mathcal H$ which corresponds to the summand 
$\sigma_i$. Then each $\mathcal H_i$ is actually 
invariant under $\sigma(K)$.

\end{lemma}
\begin{proof}
Let $a\in K$. 
Then the linear operator 
$\mathcal T$ defined as
\begin{eqnarray*}
\mathcal T:\mathcal H_1&
\mapsto& \mathcal H\\
\mathcal T(w)&=&\sigma(a)w \textrm{\ \ for all\ }
w\in\mathcal H_1
\end{eqnarray*}
is an element of 
$\mathrm{Hom}_{K'}(\sigma_1^a,\sigma)$. 
The orthogonal projection $\mathcal P$ 
of $\mathcal H$ onto $\mathcal H_2$ is an element
of $\mathrm{Hom}_{K'}(\sigma,\sigma_2)$. 
By the hypothesis of the lemma, one should have
$\mathcal P\circ\mathcal T=0$. Therefore
$\sigma(a)\mathcal H_1\subseteq \mathcal H_1$.

\end{proof}

\subsection{Statement of the main theorems}
\label{statementofthe}
We state the main theorems in this section and 
prove them in sections 
\ref{proofofthemaintheorems}
and \ref{proofofthemaintheoremsII}.\\

Throughout this section, we assume
$G,\mathbf G_\mathbb F$ and 
$\germ g$ are as in section \ref{notation} 
and $\germ g$ satisfies condition
(\ref{gdefine}).
Moreover, we assume $N_\gg$ is the
H-tower subgroup
of $G$ 
(see equality (\ref{levidecompositionofpgamma})), 
described as 
in (\ref{OKPDecomp}).\\

\begin{theorem}\label{rankablesupport}
Let $\pi$ be a unitary representation of $G$.
Then in the 
direct integral decomposition 
\[\pi_{|N_\gg}=\int_{\hat{N}_\gg}^\oplus\tau d\mu(\tau)\]
the support of the spectral 
measure $\mu$ is inside the subset
of rankable representations of $N_\gg$.
\end{theorem}

In the next theorem the subset of rankable
representations of rank $k$ of $N_\gg$ 
will be 
denoted by $\hat N_\gg(k)$.

\begin{theorem}\label{main}
Let $\pi$ be a unitary representation of $G$ 
on a Hilbert
space ${\cal H}_\pi$.
Consider the direct integral 
decomposition
\[\pi_{|N_\gg}=
\int_{\hat N_\gg}^\oplus \tau d\mu(\tau)\]
and 
let $\mathcal P_\mu$ be the projective
measure corresponding to this decomposition.
Let $\hat N_\gg(\pi)$ be the support of $\pi$
in this direct integral decomposition. 
Set
$$
{\cal H}_\pi^j=
\mathcal P_\mu(\hat{N}_{\gg}(j)\cap \hat N_\gg(\pi))
\cdot{\cal H}_\pi
\qquad\textrm{for any }0\leq j\leq \hgt(N_\gg).
$$
Then for every $j$ such that
$$
j\in\{0,1,2,...,\hgt(N_\gg)\},$$
${\cal H}_\pi^j$ is a $G$-invariant 
subspace of $\mathcal H_\pi$. 
The direct sum of all these subspaces
is equal to 
${\cal H}_\pi$.
\end{theorem}

\noindent{\bf Remark.} 
It is clear that for 
$\mathcal P_\mu(\hat{N}_{\gg}(j)\cap \hat N_\gg(\pi))$
to make sense we need to show that 
the sets $\hat N_\gg(j)
\cap \hat N_\gg(\pi)$ are indeed 
Borel subsets of $\hat N_\gg$. 
In fact one can show that the set of rankable representations
of a given rank of an H-tower group $U=H^1\cdots H^t$ 
can be constructed with a finite number of set-theoretic operations
on open and closed sets of the unitary dual of $U$.
Here is a sketch of
the proof. Any rankable representation $\sigma$ is a tensor product of 
the form 
$$
\rho_1\otimes\cdots\otimes\rho_t
$$
such that $\rho_1,...,\rho_{\mathrm{rank}(\sigma)}$ are
extensions of representations (with nontrivial central
characters)
of $H^1,...,H^{\mathrm{rank}(\sigma)}$
respectively, and the rest of the $\rho_j$'s are trivial. 
The first requirement imposes open conditions on the subset
of rankable representations of a given rank, whereas the second
requirement imposes a closed condition.

\begin{definition}\label{definitionofrank}
Let $\pi$ be any 
unitary representation of $G$. 
$\pi$ is called ``pure-rank'' if its 
restriction to $N_\gg$
is a direct integral of rankable 
representations of a fixed rank.
The common rank of these rankable representations is 
called the rank of $\pi$.

\end{definition}

Although Theorem \ref{main} is slightly 
stronger than
Corollary \ref{smallispure} below, we would like
to state it
in order to clarify the analogy between
our new theory 
and the older one.
\begin{corollary}\label{smallispure}
Let $\pi$ be an irreducible representation 
of $G$. 
Then $\pi$ is pure-rank.
\end{corollary}

\subsection{Proof of Theorem \ref{rankablesupport}}
\label{proofofthemaintheorems}
In this section we prove Theorem 
\ref{rankablesupport}.
Without loss of generality we can assume
that condition 
(\ref{ASSUMPTION}) holds.
Theorem \ref{rankablesupport} is proved
by induction on 
$\hgt(N_\gg)$.
Let $N_\gg$ be 
described as in
(\ref{OKPDecomp}).
By Howe-Moore's 
theorem \cite{HoweMoore}, when
$\hgt(N_\gg)=1$ 
there is nothing to prove.
Now let $\hgt(N_\gg)>1$. 
Let $P,L, M$ and 
$N$ be as in section
\ref{mackeyanalysissection}.
Let $P_1, M_1$ be as in Proposition \ref{definingp1}.
The H-tower 
subgroup
which is associated to $M$ 
is
$N_2=N_\gg/N$ and it
lies within $M_1$. 
But $\hgt(N_2)=\hgt(N_\gg)-1$, 
and therefore Theorem \ref{rankablesupport}
holds for $M_1$.\\

Without loss of generality, we can 
assume $\pi$ has no 
nonzero $G$-fixed vectors.
Let $\sigma_i$'s be the representations
which appear in the decomposition of
$\pi$ using
Lemma 
\ref{mackeytheoryforP}.
By a straightforward application of 
Mackey's subgroup theorem
(see Theorem \ref{MSUBT}),
or even more directly (with only little difficulty)
by writing the 
definition of the induced representation 
$\Ind_J^P\sigma_i$
explicitly, 
one can see that
the representation
$$
\Res_{P_1}^P\Ind_J^P\sigma_i
$$
is supported on representations of the form
$$
\Res_{P_1}^J\sigma_i^x \qquad x\in P
$$
where $\sigma_i^x$ is defined as in
(\ref{Repaction}); i.e.
\begin{equation}
\label{CONJUG1}
\textrm{for any }y\in J,\ \sigma_i^x(y)=
\sigma_i(x^{-1}yx).
\end{equation}
Note that $J$ is a normal subgroup of $P$
because it contains $[P,P]$. \\

By 
(\ref{sigmaidecomp}),
$\Res_{P_1}^J\sigma_i^x$ is unitarily equivalent to
$\nu_i^x\otimes \rho_i^x$ where $\nu_i^x$
and $\rho_i^x$
are defined similar to (\ref{CONJUG1}).
Theorem \ref{rankablesupport} follows from
Lemma \ref{rankonedirect} applied to 
$\rho_i^x$, the fact that
by
induction hypothesis $\nu_i^x$ is supported on 
rankable representations of $N_2=N_\gg/N$, and 
Definition \ref{s-OKPdef}.

\subsection{Proof of Theorem \ref{main}}
\label{proofofthemaintheoremsII}

We will now proceed towards
proving  Theorem \ref{main}.
Without
loss of generality, we can assume that 
condition
(\ref{ASSUMPTION}) holds.
Let $P, L, N$ and $M$ be as in section 
\ref{mackeyanalysissection}.
Let $P_1, M_1$ be as in
Proposition \ref{definingp1}.
The main idea behind the proof is that 
the Kirillov orbits of 
representations of 
different rank have
different dimensions. We will apply
some basic Kirillov theory.\\

Recall from section \ref{statementofthe}
that
$G,\germ g,\germ g_\mathbb F$ are as in 
section \ref{notation} and $\germ g$ satisfies condition
(\ref{gdefine}). Let $\Sigma,\Sigma^+$ and $\Sigma_B$ be
as in section \ref{notation} too.
Let $\tilde\beta$
be the highest root of $\germ g$ 
(see section \ref{notation}) and
$\beta$ be a simple restricted 
root such that $(\beta,\tilde\beta)=1$. 
Recall that 
\begin{equation}
\label{pbetta11}
P_{\{\beta\}}=L_{\{\beta\}}N_{\{\beta\}}
\end{equation} 
is
a standard parabolic subgroup of $G$ associated 
to $\{\beta\}$ (see section \ref{notation}).
Let $N_\gg$ be the H-tower subgroup of 
$G$ and suppose $\germ n_\gg$ is its Lie algebra. 
Let $N_{\{\beta\}}$ be as in equation
(\ref{pbetta11}) and suppose 
$\germ n_{\{\beta\}}$ 
is its Lie algebra. We consider $\germ n_\gg$ and
$\germ n_{\{\beta\}}$ 
as subalgebras of $\germ g_\mathbb F$.
Recall that for any $\gamma\in\Sigma$,
${(\germ g_\mathbb F)}_\gamma$ denotes the restricted
root subspace of $\germ g_\mathbb F$ 
associated to $\gamma$. 
Our next aim is to define a subgroup
$N_\gg^\beta$ of $N_\gg$ in case 
$\hgt(N_\gg)>1$. \\

Let $L_\gg$ be as in equation (\ref{levidecompositionofpgamma}),
with Lie algebra $\germ l_\gg$ where $\germ l_\gg\subset
\germ g_\mathbb F$, and 
$$\Sigma^L=\{\gamma\in \Sigma\ |\ {(\germ g_\mathbb F)}_\gamma
\subset \germ l_\gg\}.$$
Define
\begin{equation}
\label{newngammabeta}
\germ n_\gg^\beta=\bigoplus_{\gamma\in S_\gg} 
{(\germ g_\mathbb F)}_\gamma
\end{equation}
where 
\begin{equation}
\label{SGG}
S_\gg=\{\gamma\in\Sigma^+\ |\ {(\germ g_{\mathbb F})}_\gamma\subset 
\germ n_\gg\cap\germ n_{\{\beta\}}
\textrm{ and }\gamma-\beta\notin\Sigma^L\}.
\end{equation}
One can see that 
$\germ n_\gg^\beta$ is a Lie subalgebra of $\germ g_\mathbb F$.
To see this, assume $\gamma,\gamma'\in S_\gg$
and $\gamma+\gamma'\in\Sigma^+$. Then since both
$\germ n_\gg$ and $\germ n_{\{\beta\}}$ are Lie algebras, 
${(\germ g_\mathbb F)}_{\gamma+\gamma'}\subset \germ n_\gg\cap 
\germ n_{\{\beta\}}$. 
Moreover, if $\gamma+\gamma'-\beta\in\Sigma^L$, then 
either $\gamma\in\Sigma^L$ or $\gamma'\in\Sigma^L$
which is a contradiction. 
Consequently $\gamma+\gamma'\in S_\gg$.\\

The group $N_\gg^\beta$ is defined as the subgroup of 
$N_\gg\cap N_{\{\beta\}}$ with Lie algebra 
$\germ n_\gg^\beta$. (It is worth mentioning that
$N_\gg^\beta$ is a proper subgroup of $N_\gg\cap N_{\{\beta\}}$
if and only if $\gg$ contains an element which is not orthogonal to
$\beta$ in $\Sigma$.)\\

\begin{lemma}
\label{normalizergroupMB}
Let $\hgt(N_\gg)>1$ and $L_{\{\beta\}}$ be as in equation (\ref{pbetta11}). 
Let $M_{\{\beta\}}=[L_{\{\beta\}},L_{\{\beta\}}]$,
the commutator subgroup of $L_{\{\beta\}}$. 
Then $M_{\{\beta\}}$ normalizes $N_\gg^\beta$.
\end{lemma}
\begin{proof}
It suffices to show that for any $\gamma\in S_\gg$, if 
$\gamma+\beta\in\Sigma$ (respectively $\gamma-\beta
\in\Sigma$) then $\gamma+\beta\in S_\gg$ (respectively
$\gamma-\beta\in S_\gg$). Note that $\gamma-\beta$ cannot
be zero.\\

Assume $\gamma\in S_\gg$ and $\gamma+\beta\in\Sigma$. 
Since ${(\germ g_{\mathbb F})}_\gamma\subset \germ n_{\{\beta\}}$,
we have ${(\germ g_\mathbb F)}_{\gamma+\beta}\subset\germ n_{\{\beta\}}$
too. Similarly, since
${(\germ g_{\mathbb F})}_\gamma\subset \germ n_\gg$, 
$\gamma+\beta\notin\Sigma^L$ and
$(\gamma+\beta)-\beta=\gamma\notin\Sigma^L$ 
which imply that $\gamma+\beta\in S_\gg$.\\

Next assume $\gamma\in S_\gg$ and $\gamma-\beta\in\Sigma$.
Note that $\gamma-\beta\in\Sigma$
implies $\gamma-\beta\in\Sigma^+$
since $\gamma\in\Sigma^+$.
Again since  ${(\germ g_{\mathbb F})}_\gamma\subset \germ n_{\{\beta\}}$,
we have ${(\germ g_\mathbb F)}_{\gamma-\beta}
\in\germ n_{\{\beta\}}$. Moreover, by the
definition of $\germ n_\gg^\beta$, $\gamma-\beta\notin\Sigma^L$
and therefore ${(\germ g_{\mathbb F})}_{\gamma-\beta}
\subset\germ n_\gg$. Finally, if 
$(\gamma-\beta)-\beta\in\Sigma^L$,
then $\gamma=2\beta+\gamma_1$ where $\gamma_1\in\Sigma^L$. This means 
that $(\gamma,\tilde\beta)=2$, i.e. $\gamma$ is the highest root.
Since any simple restricted root appears in the highest root with a 
positive coefficient, it follows that $\Sigma_B$ 
consists of $\beta$ and the elements of $\gg$.
Consequently, there is only one simple restricted root outside $\gg$, which
implies that $\hgt(N_\gg)=1$, which contradicts our assumption.

\end{proof}

\begin{proposition}\label{ordim2}
Let $G$ be as in Theoreom \ref{main},
$N_\gg$ be the H-tower subgroup of $G$ and
$\hgt(N_\gg)>1$.
Let $N_\gg^\beta$ be defined as 
in equation (\ref{newngammabeta}).
Then the 
restriction of a rankable representation 
$\sigma$ of rank $k$ 
of $N_\gg$ to
$N_\gg^\beta$
is a direct integral of irreducible 
representations of $N_\gg^\beta$
whose attached coadjoint orbits have 
the same dimension equal to $2(n_1+\cdots +n_k-c)$, 
where
$c$ is the codimension of $N_\gg^\beta$ in $N_\gg$, 
and $n_i$'s are defined as in
Corollary \ref{OKPORBITSDIM}.

\end{proposition}

\begin{proof}
Throughout the proof we assume 
$\mathbb F=\mathbb R$ for
simplicity. The proof for other 
local fields is essentially 
similar.\\
  
We will analyze
$$\mathrm{Res}_{N_\gg^\beta}^{N_\gg}(\rho_{1}
\otimes\cdots\otimes\rho_{k})$$
with the $\rho_{i}$'s as in 
Definition \ref{s-OKPdef}.
It is easy to see that
$N_\gg^\beta\supset H^2\cdots H^{\hgt(N_\gg)}$,
so when 
$j>1$ the restriction of $\rho_{j} $
to ${N_\gg^\beta}$ 
is irreducible. It remains to understand the 
restriction of 
$\rho_{1}$ to $N_\gg^\beta$. 
The group $H^1\cap N_\gg^\beta$ 
is a direct product of 
a Heisenberg 
group of dimension $2(n_1-c)+1$ 
and a $c$-dimensional 
abelian 
group whose Lie algebra corresponds 
to the direct sum of restricted root spaces
$(\germ g_\mathbb R)_{\tilde{\beta}-\gamma}$, 
such that 
$(\gamma,\tilde\beta)>0$ 
but $\gamma\notin S_\gg$.
To see why these restricted root spaces 
form an isotropic subspace of the Heisenberg
nilradical of $\germ g_\mathbb R$, 
suppose $\tilde\beta-\gamma_1,\tilde\beta-\gamma_2$
are given such that for any $i\in\{1,2\}$,
$(\gamma_i,\tilde\beta)>0$ 
but $\gamma_i\notin S_\gg$.
If $(\tilde\beta-\gamma_1)+(\tilde\beta-\gamma_2)=\tilde\beta$,
then $\tilde\beta=\gamma_1+\gamma_2$. But for any $i\in\{1,2\}$, 
$\gamma_i=\beta+\gamma_i'$ where $\gamma_i'\in \Sigma^L$
or $\gamma_i'=0$.
Therefore $\tilde\beta=2\beta+\gamma_1'+\gamma_2'$, which
implies that $\Sigma_B$ consists of $\beta$ and the elements of $\gg$.
(See the proof of Lemma \ref{normalizergroupMB}.) Consequently
$\hgt(N_\gg)=1$, which is a contradiction.\\


Let the 
decomposition of
$H^1\cap N_\gg^\beta$ as a direct product be 
$H^\beta\times \mathbb R^c$,
where $H^\beta$ is the
$2(n_1-c)+1$-dimensional Heisenberg group
and $\mathbb R^c$ is the abelian factor. 
We will denote the irreducible
representation of $H^\beta$
with central character $\chi_1$ by
$\rho^\beta_1$.
\begin{lemma}\label{rho1} 
Under the foregoing assumptions,
\[\mathrm{Res}^{N_\gg}_{H^1\cap N_\gg^\beta}
\rho_{1}=\int_{{\mathbb R^c}^*} 
\rho^\beta_{1}\hat{\otimes}
\psi_s 
d\mu(s)\]
where each $\psi_s(t)$, given by 
$\psi_s(t)=e^{\mathrm is(t)}$
for some $s\in {\mathbb R^c}^*$,
the vector-space dual of $\mathbb R^c$,
is a unitary character 
of $\mathbb R^c$. 
\end{lemma}

\begin{proof}
Clearly
$$\mathrm{Res}^{N_\gg}_{H^1\cap N_\gg^\beta}
\rho_{1}=
\mathrm{Res}_{H^\beta\times 
\mathbb R^c}^{H^1}
\mathrm{Res}^{N_\gg}_{H^1}\rho_{1}=
\mathrm{Res}_{H^\beta
\times \mathbb R^c}^{H^1}
\rho_{1}.$$
The space ${\cal H}_{\rho}$ 
of any representation
$\rho$ of the Heisenberg group 
$H^1=H_{n_1}$
introduced in section 2.2
can be written 
as 
\begin{eqnarray*}
{\cal H}_{\rho}=
L^2(\germ y)=
L^2(\mathrm{Span}_\mathbb R\{Y_1,\ldots,Y_c\}
\oplus   
\mathrm{Span}_\mathbb R\{Y_{c+1},\ldots,Y_{n_1}\})=\\
L^2(\mathbb R^c)
\hat{\otimes}L^2(\mathrm{Span}_\mathbb R
\{Y_{c+1},\ldots,Y_{n_1}\})
\approx
\int_{{\mathbb R^c}^*}^\oplus L^2_s d\mu(s)
\end{eqnarray*}
where each
$L^2_s$ is equal to 
$L^2(\mathrm{Span}_\mathbb R\{Y_{c+1},
\ldots,Y_{n_1}\})$ 
on which $\mathbb R^c$ 
acts via the character $\chi_s(x)=e^{\mathrm is(x)}$.
In fact $L^2_s$ is a representation of
$H^\beta\times\mathbb R^c$.
Lemma \ref{rho1} is proved.
\end{proof}

Next we show that $\mathbb R^c\subset
\mathcal Z(N_\gg^\beta)$. To see this, take a restricted
root space ${(\germ g_\mathbb R)}_{\tilde\beta-\gamma}$
which lies inside $\mathbb R^c$. It suffices to show that for
any $\gamma'\in S_\gg$, $(\tilde\beta-\gamma)+\gamma'\notin\Sigma$.
But we know that $\gamma=\beta+\gamma_1$ where
$\gamma_1\in\Sigma^L$ or $\gamma_1=0$. 
Therefore $(\tilde\beta-\gamma)+\gamma'=
\tilde\beta-\beta-\gamma_1+\gamma'$. Consequently, if
$(\tilde\beta-\gamma)+\gamma'\in\Sigma$, 
then either $\gamma'\in\Sigma^L$
or $\gamma'-\beta=\gamma_1\in\Sigma^L$. However, none of these 
is possible for $\gamma'$
by the definition of $S_\gg$ in (\ref{SGG}).\\

We have shown that $\mathbb R^c\subset
\mathcal Z(N_\gg^\beta)$
and
the action of $\mathbb R^c$ on 
each $L_s^2$ is by a distinct character $\psi_s$
(see Lemma \ref{rho1}).
Hence the restriction of ${\rho_{1}}$ to ${N_\gg^\beta}$ 
breaks into a direct
integral of a $c$-parameter family of 
irreducible representations. Consequently
the same thing happens 
to any rankable representation 
$$\sigma=\rho_{1}
\otimes\cdots\otimes\rho_{k}.$$
By Theorem 2.5.1 of \cite{CG}
and Lemma
\ref{N_1N_2} applied to
$N_\gg^\beta$ which is a semidirect
product of $H^\beta\times \mathbb R^c$ 
and
$N_2$,
it follows that the projection of 
the
coadjoint orbit $\mathcal O^*_\sigma$ onto 
the Lie algebra of 
$N_\gg^\beta$ is foliated by subvarieties of 
codimension $2c$, which are 
indeed coadjoint orbits of the constituents
of the rankable representations in the  
restriction of $\sigma$ to $N_\gg^\beta$.
Proposition \ref{ordim2} is proved.
\end{proof}

We will now concentrate on finishing the proof
of Theorem \ref{main}.
Theorem \ref{main} is proved by induction on
the height $\hgt(N_\gg)$
of the 
H-tower group $N_\gg$.
If $\hgt(N_\gg)=1$ then there is nothing 
to prove. Let 
$\hgt(N_\gg)>1$. Then the H-tower
subgroup of $M_1$ 
is
equal to $N_2$, where
$$
N_2=N_\gg/N.
$$ 
Clearly $\hgt(N_2)=\hgt(N_\gg)-1$.\\

Without loss of generality we can assume $\pi$
has no nonzero $G$-invariant vectors. 
Consider the decomposition of $\pi$ given in
(\ref{mackeydecompositionofP}).
Applying induction hypothesis to $M_1$,
which contains the H-tower subgroup
$N_2$,
we can refine this decomposition by
expressing each $\nu_i$ as a direct sum of
its $M_1$-invariant 
pure-rank parts (where the rank 
for a representation of $M_1$ 
is naturally defined
with respect to $N_2$). 
Conseqently, as a representation of $M_1$,
$$
\nu_i=\nu_{i,0}\oplus\nu_{i,1}
\oplus\cdots\oplus
\nu_{i,\hgt(N_2)}.
$$
Here $\nu_{i,j}$ denotes the part of $\nu_i$
supported on rankable representations of 
$N_2$ of 
rank $j$. (Note 
that some of the $\nu_{i,j}$'s may be
trivial, but this fact does not affect our proof.)
Therefore we have
$$
{\sigma_i}_{|P_1}=(\nu_{i,0}\otimes
\rho_i)\oplus\cdots\oplus(\nu_{i,\hgt(N_2)}
\otimes\rho_i).
$$
Let $\mathcal H$ be the Hilbert space of
the representation $\sigma_i$.
For any  $j\in\{0,1,2,...,\hgt(N_2)\}$,
let  $\mathcal H_{j}$ be the subspace of
$\mathcal H$ which corresponds to
$\nu_{i,j}\otimes \rho_i$. 
Our next task is to prove that each
$\mathcal H_j$ is in fact invariant under
$\sigma_i(J)$. 
(See Lemma \ref{mackeytheoryforP}.)
This follows from Lemma
\ref{orbitinvariance} once we prove the following lemma.

\begin{lemma}
\label{LEMMAforseparation}

Let $\eta=\nu_{i,j}\otimes\rho_i$ and 
$\eta'=\nu_{i,j'}\otimes\rho_i$ where 
$j\neq j'$. Let $a\in J$. Then
$$\mathrm{Hom}_{P_1}(\eta^a,\eta')=\{0\}$$
where $\eta^a$ is a representation of 
$P_1$ defined on $\mathcal H_j$ defined as
in (\ref{Repaction}); i.e.
$$\eta^a(x)=\eta(a^{-1}xa)\quad\textrm{for all\ }
x\in P_1.
$$
\end{lemma}

\begin{proof}
We actually prove more; i.e. that 
\begin{equation}
\label{homforngamma}
\mathrm{Hom}_{N_\gg}(\Res^{P_1}_{N_\gg}\eta^a,
\Res_{N_\gg}^{P_1}\eta')
=\{0\}.
\end{equation}
We claim that
$\Res_{N_\gg}^{P_1}\eta^a$ is a direct integral
supported on 
rankable representations of $N_\gg$
of rank $j+1$. The claim implies
(\ref{homforngamma})
because if $j\neq j'$, then it implies that
$\Res^{P_1}_{N_\gg}\eta^a$ and 
$\Res^{P_1}_{N_\gg}\eta'$
are supported on disjoint subsets of 
$\hat N_\gg$. We prove this claim below.\\

By Lemma \ref{rankonedirect}, 
$\rho_i^a$ is supported on
rankable representations of $N_\gg$ of rank one.
Next we show that $\nu_{i,j}^a$ is supported
on rankable representations of $N_2=N_\gg/N$ of
rank $j$.\\

As we know, $G$ is a central extension of
$\mathbf G_\mathbb F$. Suppose this extension 
is represented as in (\ref{centralextensionofG}). \\

First let $a\in N$. Then $a\in N_\gg$ and consequently
$\Res_{N_\gg}^{P_1}\eta^a$ is unitarily equivalent to
$\Res_{N_\gg}^{P_1}\eta$. Since $J\subset P$ and any
element of $P$ is a product of an element of 
$N$ and 
an element of $p^{-1}(\mathbf L_\mathbb F)$, 
it follows that it
suffices to assume
that 
$p(a)\in \mathbf L_\mathbb F$.\\

It follows from $p(a)\in
\mathbf L_\mathbb F$ that   
$p(a)\mathbf M_\mathbb F p(a^{-1})=
\mathbf M_\mathbb F$. The group
$p(P_\gg)$ is the 
$\mathbb F$-points of the 
$\mathbb F$-parabolic $\mathbf P_{\mathbf\gg'}$
of $\mathbf G$. Now
$\mathbf P_{\mathbf m}=
\mathbf P_{\mathbf\gg'}\cap \mathbf M$
is an $\mathbb F$-parabolic of $\mathbf M$.
Since $p(a)\in\mathbf L_\mathbb F$,
$p(a)\mathbf P_\mathbf m p(a^{-1})$ is another
$\mathbb F$-parabolic of 
$\mathbf M$, and therefore
by \cite[Theorem 20.9]{Borel} these two parabolics
are conjugate under an element $p(b)$ of 
$\mathbf M_\mathbb F$ where 
$b\in M_1$ (recall 
from Proposition \ref{definingp1}
that 
$p(M_1)\supseteq \mathbf M_\mathbb F$). 
This means that 
$p(ba)\mathbf P_\mathbf m p(ba)^{-1}=
\mathbf P_\mathbf m$.
If $\mathbf U_\mathbf m$ is the 
unipotent radical of
$\mathbf P_\mathbf m$, then 
$p(ba)\mathbf U_\mathbf m p(ba)^{-1}=
\mathbf U_\mathbf m$, which 
implies that $baN_2a^{-1}b^{-1}=N_2$.\\

Let $c=ba$. Then $c\in P$.
Let $\germ n_2\subset \germ n_\gg$ 
be the Lie algebra of $N_2$.
Consider $Ad^*(c)$ as 
a linear map from the dual of the Lie
algebra of $P$ to itself.
Obviously $Ad^*(c)(\germ n_2^*)=\germ n_2^*$.
Let $\tau$ be an irreducible 
unitary representation of $N_2$ and let 
$\mathcal O^*_\tau$ be the coadjoint orbit
associated to $\tau$. Then
the coadjoint
orbit associated to $\tau^c$ is 
$Ad^*(c)(\mathcal O^*_\tau)$. This fact follows 
for instance from 
\cite[III 11]{dufloacta}. A short proof of
this fact can be obtained by an adaptation of
the proof of Lemma \ref{lemmaondimsagain}.\\


For simplicity let 
$\nu=\nu_{i,j}$ for fixed $i,j$.
The representation $\Res_{N_2}^{M_1}\nu$
is supported on rankable representations 
of $N_2$ whose associated coadjoint orbits have 
dimension 
$$n_2+\cdots+n_{j+1}.$$ 
(See Corollary
\ref{OKPORBITSDIM}
for the definition of $n_i$'s.)\\
 
Let $\nu^a,\nu^b,\nu^c$ be representations of
$M_1$
defined as in (\ref{Repaction}).
(Note that $aM_1a^{-1}=M_1$ because 
$P_1$ is normal in 
$p^{-1}(\mathbf P_\mathbb F)$.)
Since $b\in M_1$, any representation $\theta$ 
of $M_1$ is obviously unitarily equivalent
to $\theta^b$. 
Therefore $\nu^c$ is unitarily 
equivalent to $\nu^a$.
Since the action of
$Ad^*(c)$ on $\germ n_2^*$ is linear, it does not change
the dimension of coadjoint orbits, and therefore 
$\Res^{M_1}_{N_2}\nu^c$ is supported on  
unitary representations of $N_2$ whose associated
coadjoint
orbits have dimension $n_2+\cdots+n_{j+1}$.\\

Since $\nu^c$ is a unitary representation 
of $M_1$, by Theorem \ref{rankablesupport} 
all representations in
the support of its restriction to $N_2$ should be
rankable. Therefore $\Res^{M_1}_{N_2}\nu^c$ is 
supported on rankable representations of $N_2$ 
of rank $j$. Now 
$\nu^a$ is unitarily equivalent to $\nu^c$
as a representation of $M_1$, and hence as a 
representation of $N_2$. Therefore $\nu^a$
is also supported on rankable representations of 
$N_2$ of rank $j$. Definition \ref{s-OKPdef} completes
the proof of our claim. 
The proof of
Lemma \ref{LEMMAforseparation} is complete.

\end{proof}

As mentioned before, Lemma \ref{LEMMAforseparation}
implies that
each of the 
subspaces $\mathcal H_j$ (which corresponds
to the representation $\nu_{i,j}\otimes\rho_i$ of
$P_1$)  
is invariant under
$\sigma_i(J)$. 

\begin{lemma}
\label{pionP}
There is a $P$-invariant direct sum
decomposition 
of $\pi$ such as
\begin{equation*}
\pi_{|P}=\pi_1\oplus\cdots\oplus\pi_{\hgt({N_\gg})}
\end{equation*}
where each $\pi_i$ is of pure $N_\gg$-rank $i$.\\
\end{lemma}

\begin{proof}
The invariance of the spaces $\mathcal H_j$ under
$\sigma_i(J)$
means that one
can actually decompose $\sigma_i$ as a direct
sum 
$$
\sigma_i=\sigma_{i,1}\oplus\sigma_{i,2}\oplus\cdots
\oplus\sigma_{i,\hgt(N_\gg)}
$$
such that each $\sigma_{i,j}$ is a representation of $J$
and 
$$
{\sigma_{i,j}}_{|P_1}=\nu_{i,j-1}\otimes \rho_i.
$$
Now
an application of Mackey's subgroup 
theorem (see Theorem \ref{MSUBT})
immediately implies that  
the representation
$$
\Res_{N_\gg}^P
\mathrm{Ind}_{J}^P\sigma_{i,j}
$$
is supported on rankable representations
of $N_\gg$ of rank $j$. (See the argument of Theorem
\ref{rankablesupport} and Lemma \ref{LEMMAforseparation}.)
Lemma \ref{pionP} is proved.
\end{proof}


Recall from Lemma
\ref{normalizergroupMB} that
the group $N_\gg^\beta$ is normalized 
    by $M_{\{\beta\}}$. 
A version of the following lemma was mentioned
in the proof of Lemma \ref{LEMMAforseparation}.

\begin{lemma}
\label{lemmaondimsagain}
Let $\tau$ be a unitary representation of 
$N_\gg^\beta$ and let $g\in M_{\{\beta\}}$. 
Let $\tau^g$ be defined as in 
(\ref{Repaction}).
Then
the coadjoint orbits associated to $\tau$ and $\tau^g$
(in the sense of Kirillov's orbital theory)
have the same dimension.
\end{lemma}

\begin{proof}
The proof follows immediately from the fact that
the coadjoint orbit attached to $\tau^g$
is $Ad^*(g)(\mathcal{O}^*)$, which follows from
\cite[III 11]{dufloacta}. 
We would like to give a short proof of this fact 
for the reader's convenience. 
Let $\germ n_\gg^\beta$ be the Lie algebra
of $N_\gg^\beta$. 
Fix an additive character $\chi$ of 
    $\mathbb F$ as done in \cite[\S4]{moore}. 
(When $\mathbb F=\mathbb R$ or $\mathbb C$,
    $\chi(t)=e^{\mathrm{i}\mathrm{Re}(t)}$ 
where $\mathrm{Re}(t)$ means the real part of $t$,
and when $\mathbb F$ is $p$-adic, $\chi$
will be an unramified character given by Tate.)
Let the coajoint orbit associated to $\tau$ be
$\mathcal O^*_\tau\subset \germ n_\gg^\beta$. 
By Kirillov's orbital theory we know that
$\tau$ is constructed as follows.
    One chooses an arbitrary element 
    $\lambda\in \mathcal O^*_\tau$ 
and a maximal 
    subalgebra $\germ q$ of $\germ n_\gg^\beta$
        subordinate to $\lambda$,
    which exponentiates to a closed 
	subgroup $Q$ of $N_\gg^\beta$.  
    Then
    $$\tau=\mathrm{Ind}_{Q}^{N_\gg^\beta}
(\chi\circ{\lambda\circ\log}).$$
    For any $g\in M_{\{\beta\}}$, let 
$(\chi\circ{\lambda\circ\log})^g$
be defined as in (\ref{Repaction}).
Then
    \begin{eqnarray*}
    \tau^g=
    \mathrm{Ind}_{gQg^{-1}}^{N_\gg^\beta}
    (\chi\circ{\lambda\circ\log})^g=
    \mathrm{Ind}_{gQg^{-1}}^{N_\gg^\beta}
    (\chi\circ{Ad^*(g)(\lambda)\circ\log}).
    \end{eqnarray*}
    Since $Ad(g)(\germ q)$ is a
    maximal subalgebra subordinate to
    $Ad^*(g)(\lambda)$,
    the coadjoint
    orbit attached to $\tau^g$  
    is $$Ad^*(g)(\mathcal O^*).$$ 
    Since the action of $Ad^*(g)$ is linear on 
    $\germ n_\gg^\beta$, it does not change the dimension
of the coadjoint orbit.
\end{proof}

To prove Theorem \ref{main}, 
we show 
that each of the components $\pi_i$ 
given in Lemma \ref{pionP}
is $G$-invariant as well. To this
end, we first 
prove the following lemma.

\begin{lemma}
\label{lemlemlemlem}
The $P$-invariant 
decomposition in
Lemma \ref{pionP} is preserved by 
the action of $M_{\{\beta\}}$, where 
$$M_{\{\beta\}}=[L_{\{\beta\}},L_{\{\beta\}}]$$ 
and $L_{\{\beta\}}$ is the Levi component 
of the standard parabolic 
$P_{\{\beta\}}$ with $\beta$ as in
Proposition \ref{ordim2}.

\end{lemma}

\begin{proof}
    
    For any $1\leq j\leq \hgt(N_\gg)$, 
the representation
    $$\Res^P_{N_\gg^\beta}\ \pi_j$$
    is a direct integral of
    representations which correspond to coadjoint orbits
    of dimension  
    $2(n_1+\cdots+n_j-c)$. 
    However, $n_1,...,n_{\hgt(N_\gg)}>0$.
Therefore, if we define 
$\pi_j^a$ (for any $a\in M_{\{\beta\}}$)
as in (\ref{Repaction}), then
Lemma \ref{lemmaondimsagain}
implies that
for any $j'\neq j$, the dimension of the
coadjoint orbits associated to the 
irreducible representations 
of $N_\gg^\beta$
in the support of $\Res_{N_\gg}^P\pi_j^a$ is different from the 
dimension of the coadjoint orbits assoiated to
the irreducible 
representations of $N_\gg^\beta$
in the support of $\Res_{N_\gg}^P\pi_{j'}$. This means that
$$
\mathrm{Hom}_{N_\gg^\beta}(\Res_{N_\gg^\beta}^P \pi_j^a,
\Res_{N_\gg^\beta}^P\pi_{j'})=\{0\}.
$$
Now we apply 
Lemma \ref{orbitinvariance}
with $K=[P_{\{\beta\}},P_{\{\beta\}}],
K'=N_\gg^\beta, \sigma=\pi,
\sigma_1=\pi_i$ and 
$$\sigma_2=\bigoplus_{i\neq j}\pi_i.$$
It follows that each $\pi_j$ is 
invariant under the action of $M_{\{\beta\}}$.
Lemma \ref{lemlemlemlem} is proved.

\end{proof}

 To finish the proof of Theorem \ref{main}, 
we note that the
 parabolic subgroup $P$ in $G$ 
 is maximal, therefore the 
 group 
 generated by $M_{\{\beta\}}$ and $P$ will be 
 equal to $G$. (This follows from the 
Bruhat-Tits decomposition.)
The decomposition of $\pi$ given in
Lemma \ref{pionP} is preserved by both $P$ and 
$M_{\{\beta\}}$, and hence by the group 
generated by them.
Therefore the decomposition of Lemma
\ref{pionP} is $G$-invariant.\\

\noindent{\bf Remark.\ } 
Let $G$ be as in Theorem \ref{rankablesupport}, 
such that 
the central extension identifying $G$ is
as in (\ref{centralextensionofG}).
Let $N_\gg$ be the
H-tower subgroup of the group $G$. 
Let $\pi$ be 
a unitary representation of $G$ of pure rank $k$,
where $k\leq \hgt(N_\gg)$.
Consider the subgroup 
$A_k$ of $G$ defined as follows.
$A_k=p^{-1}(\mathbf A'_\mathbb F)$, where $\mathbf A'$
is an $\mathbb F$-torus inside the maximal split 
$\mathbb F$-torus $\mathbf A$ whose Lie algebra is 
spanned by the coroots 
$H_{\tilde\beta_1},...,H_{\tilde\beta_k}$.
Let $\tau$ be a rankable representation of $N_\gg$
of rank $k$. The stabilizer $S_\tau$ of $\tau$ inside
$A_k$ is a finite subgroup of $A_k$. Moreover, 
under the action of $A_k$ there are
only a finite number of orbits of rankable 
representations of rank 
$k$. By Mackey theory, the restriction of $\pi$ to 
$A_k\ltimes N_\gg$ is a direct integral of 
representations
of the form
\begin{equation}
\label{inductionhowe}
\mathrm{Ind}_{S_\tau\ltimes N_\gg}^{A_k\ltimes N_\gg}
\sigma_\tau
\end{equation}
where $\sigma_\tau$ is irreducible and
${\sigma_\tau}_{|N_\gg}=n_\tau\tau$ for some
$n_\tau\in\{1,2,3,...,\infty\}$. 
Moreover,
by Frobenius reciprocity, 
$\sigma_\tau$ will be a subrepresentation
of $\mathrm{Ind}_{N_\gg}^{S_\tau\ltimes N_\gg}\tau$. 
Since $S_\tau$ is a finite group,
there are only a finite number of 
possibilities for $\sigma_\tau$. Since the 
number of orbits 
(and hence stabilizers) of rankable 
representations of rank $k$ 
under the action of $A_k$
is finite, the number of 
representations of the form
(\ref{inductionhowe}) is finite as well. 
This implies the following
result.
\begin{proposition}
Let $\pi$ be a unitary representation of $G$ 
of pure rank $k$.
Then there exists a finite family
$\{\tau_1,...,\tau_t\}$ of 
irreducible representations of 
$A_k\ltimes N_\gg$, independent of $\pi$, 
such that 
$$\pi_{|A_k\ltimes N_\gg}=n_1\tau_1+\cdots+n_t\tau_t$$
where $n_i\in\{1,2,3,...,\infty\}$ for each $i$.
\end{proposition}

\section{Relation with the old theory}

\subsection{Outline of the old theory}

In this section we show how the notion of rank
defined in the past
sections relates to
the existing theory 
for classical groups. 
To this end, we show that for the real forms of 
classical groups, the two notions of rank 
(the one defined
in \cite{LiInv}
and the one defined in Definition
\ref{definitionofrank}) 
are equivalent.
Here we give a brief
outline of the old theory.
In classical cases, rank of a representation of the
real semisimple group $G$ is
defined in terms of its restriction to the 
centers of nilradicals of maximal 
parabolic subgroups. One can characterize
each of these parabolics with a node in the
Dynkin diagram of the restricted root system
in a natural way. It turns out that there is a
(not necessarily unique) standard 
parabolic which 
provides the most refined information about the
rank. We will devote this section to exhibiting
the coincidence of the two notions of rank on this
parabolic. Really the main idea is some slight
modification of the 
fact that the nilradical of the rank 
parabolic subalgebra 
contains a maximal isotropic subspace of 
each of the Heisenberg algebras in 
the H-tower $N_\gg$. 
Our presentation of the results
follows the notation of older literature
\cite{Howe}, \cite{LiInv}, \cite{Scara}.\\

The notation used in this section is chosen independent of other sections
in order to simplify matters and be more coherent with older works.
For simplicity we only consider the case $\mathbb F=\mathbb R$.
The general case is essentially the same and will only be
more technical. 
It is more convenient to consider  
classical groups of different types (in the sense of
\cite{Howetheta}) separately.
Here we quickly review the definition of classical groups of 
types I and II over a local field, but later we will 
retain our assumption that $\mathbb F=\mathbb R$.\\

Let $\mathbb F$ be a local field, $D$ a division 
algebra over $\mathbb F$ with an involution,
and $V$ a left vector space over $D$ of dimension 
$n$. A classical group $G$ is said to be of type II if
$G=GL_D(V)$. From now on,
by $(\cdot,\cdot)$ we mean a 
Hermitian or 
skew-Hermitian sesquilinear 
form $(\cdot,\cdot)$ on $V$.
A classical group $G$ of type I is
the connected component of identity of the stabilizer
subgroup of $(\cdot,\cdot)$ inside $GL_D(V)$. 
The real groups of type I which are of our interest here,
i.e. those which satisfy the assumptions of 
Proposition \ref{complexheisenbergparabolicstructure},
correspond to the cases where
$\mathbb F=\mathbb R$ and $D=\mathbb{R,C}$ or 
$\mathbb H$ with their usual involutions. \\

\subsection{Groups of type I}
A typical maximal parabolic of these groups can 
be described as follows. Take a maximal 
polarization inside $V$, i.e. a maximal set 
of vectors 
$$\{e_1,\ldots,e_r,e_1^*,\ldots,e_r^*\}$$
in $V$
which satisfy 
\begin{center}   
 $(e_i,e_j)=(e_i^*,e_j^*)=0$\\
  $(e_i,e^*_j)=\delta_{i,j}.$
\end{center}
In fact $r$ is equal to the split rank of $G$.
For any $k$
let $$X_k=\mathrm{Span}_D\{e_1,\ldots,e_k\}
\quad,\quad
X_k^*=\mathrm{Span}_D\{e_1^*,\ldots,e_k^*\}$$
and define $V_k$ to be $X_k\oplus X_k^*$.
Let $P_k$ be the subgroup of $G$ that
consists of elements which
leave the subspace $X_k^*$ invariant. 
$P_k$ is a parabolic subgroup and
the Levi decomposition of $P_k$ looks like
\begin{equation}
\label{levidpk}
P_k=GL_D(X_k^*)G(V_k^\perp)N_k
\end{equation}
where by $G(V_k^\perp)$ we mean 
the stabilizer of $(\cdot,\cdot)$
as a form on $V_k^\perp$. Here $N_k$
is the unipotent radical of $P_k$.\\

$P_r$
is the parabolic which provides the most refined
information about the rank (in the 
sense of \cite{Howe},\cite{LiInv}). 
\begin{definition}
Let $r$ be the split rank 
of $G$. The parabolic $P_r$
or its Lie algebra are called the {\it rank parabolic}.\\
\end{definition}

The unipotent radical $N_k$
is a two-step nilpotent simply connected
Lie group and therefore it can be 
identified with
its Lie algebra via the exponential map. 
From now on, we think of any $N_k$ through this
identification, and although slightly 
ambiguous, we use the same notation
for its Lie subgroups and their Lie algebras.
This is done in order to avoid complicated
notation and to keep the presentation as close
to the style used in the papers of Howe and Li. 
We will make it clear whether or not we are using
a Lie group or a Lie algerba wherever necessary.\\

As in \cite{LiInv}, we have the following exact 
sequence of Lie algebras:
\begin{displaymath}0
    \longrightarrow ZN_k\longrightarrow
    N_k\longrightarrow\mathrm{Hom}_D(V_k^\perp,X_k^*)
    \longrightarrow 0
\end{displaymath}
where $ZN_k$ is the center of $N_k$.
$ZN_k$ is isomorphic to 
$$\mathrm{Hom}_D^{\mathrm{inv}}(X_k,X_k^*)$$ 
where $\mathrm{Hom}_D^{\mathrm{inv}}(X_k,X_k^*)$
is the $\mathbb F$-subspace of 
elements $T$ of $\mathrm{Hom}_D(X_k,X_k^*)$
satisfying 
$$\forall i,j\in\{1,2,\ldots,k\},
\ \ (Te_i,e_j)+(e_i,Te_j)=0.$$
Thus as an $\mathbb F-$vector space, the Lie algebra
$N_k$ 
can be expressed as
\begin{equation}
\label{lastmin}
N_k=\mathrm{Hom}_D(V_k^\perp,X_k^*)\oplus
\mathrm{Hom}_D^{\mathrm{inv}}(X_k,X_k^*).
\end{equation}

The isomorphism of $\mathrm{Hom}_D(V_k^\perp,X_k^*)$
into the Lie algebra $N_k$ can be depicted as 
\begin{eqnarray}
\label{transpose}
\tilde{}:\mathrm{Hom}_D(V_k^\perp,X_k^*)
&\mapsto&\mathrm{Hom}_D(V,V)
\end{eqnarray}
$$\ \ \ \ \ \ T\ \mapsto\ \tilde{T}$$

where $\tilde T$ is defined as follows:

\begin{eqnarray*}
\tilde{T}v=Tv\ & \textrm{for}\ & v\in V_k^\perp\\ 
\tilde{T}v=0\ \ &\textrm{for}\  & v\in X_k^*\\
\tilde{T}v=T^tv\ &\textrm{for}\ & v\in X_k.
\end{eqnarray*}
Here
$T^t\in\mathrm{Hom}_D(X_k,V_k^\perp)$ 
is defined uniquely by
\begin{eqnarray*}
    \forall\ v\in V_k^\perp\ ,\ x\in X_k
    \qquad (Tv,x)+(v,T^tx)=0.
\end{eqnarray*}

It turns out that the Heisenberg parabolic $P$ of
$G$ is $P_{k_1}$, where
$k_1=2$ for $G=SO_{p,q}$ and $k_1=1$
for all other classical cases 
under consideration. Let $P=LN$ be the Levi
decomposition of $P$. Let $M$ 
be the appropriate simple isotropic 
factor 
of $[L,L]$; i.e. we drop the 
redundant
factor of $[L,L]$ which, in
(\ref{levidpk}), corresponds to $GL_D(X_{k_1}^*)$. 
Let $\germ m$ be the Lie algebra of $M$.
For any $k$ define 
$$Y_k=\mathrm{Span}\{e_k,\ldots,e_r\}\quad,\quad
Y_k^*=\mathrm{Span}\{e_k^*,\ldots,e_r^*\}.$$
The center of the 
nilradical of the rank parabolic 
of $\germ m$ 
is identical to 
$$\mathrm{Hom}^{\mathrm{inv}}_D
(Y_{k_1+1},Y_{k_1+1}^*).$$
The Lie algebra $\mathrm{Hom}_D^{\mathrm{inv}}
(Y_{k_1+1},Y_{k_1+1}^*)$ 
acts on the Lie algebra 
$N_{k_1}$ through the adjoint action of 
$\germ m$.
By Theorem \ref{complexheisenbergparabolicstructure}
this action will be trivial on $ZN_{k_1}$. The following
simple 
lemma describes this action more explicitly.
\begin{lemma}\label{actionHOM}
Let $X\in
\mathrm{Hom}_D^{\mathrm{inv}}(Y_{k+1},Y_{k+1}^*)$.
Let $Y\in\mathrm{Hom}_D(V_k^\perp,X_k^*)$.
Then the 
adjoint action of $\germ m$
on $\germ n=N_k$ is described as 
$$\mathrm{ad}_X(\tilde Y)=-\tilde{(YX)}$$
where $\tilde{YX}$ is defined as in (\ref{transpose}).
\end{lemma}
\noindent{\bf Remark.}\ Note 
that we think of $-YX$ as an element of
$\mathrm{Hom}_D(V_k^\perp, X_k^*)$ which is zero on 
$V_r^\perp$ and $Y_{k+1}^*$.\\

The restriction of $(\cdot,\cdot)$ 
to $V_r^\perp$ is a definite form.
Without loss of generality, we may 
assume 
that the form is positive definite.
Let $\{f_1,\ldots,f_{n-2r}\}$ be an orthonormal 
basis for $V_c=V_r^\perp$; consequently 
$V_c=\mathrm{Span}_D\{f_1,\ldots,f_{n-2r}\}$.
One can see that $\mathrm{Hom}_D(V_k^\perp,X^*_k)$
is equal to
\begin{equation}\label{thisdirect}
        \mathrm{Hom}_D(Y_{k+1},X^*_k)
    \oplus
    \mathrm{Hom}_D(Y_{k+1}^*,X^*_k)
    \oplus
    \mathrm{Hom}_D(V_c,X^*_k).
\end{equation}
We consider the direct sum decomposition
(\ref{thisdirect}) inside $N_k$ (see (\ref{lastmin}) )
We observe that:
\begin{lemma}\label{commutingspaces}
In the direct sum 
decomposition (\ref{thisdirect}), the third summand
commutes with the first two summands, and the
adjoint action of 
$\mathrm{Hom}_D^{\mathrm{inv}}(Y_{k+1},Y_{k+1}^*)$
on the third summand is trivial.\\ 
\end{lemma} 
Now take $k=k_1$. Then the direct sum
\begin{equation}\label{polarization}
    \mathrm{Hom}_D(Y_{k_1+1},X^*_{k_1})
    \oplus
    \mathrm{Hom}_D(Y_{k_1+1}^*,X^*_{k_1})
    \oplus
    \mathrm{Hom}_D^{\mathrm{inv}}(X_{k_1},X_{k_1}^*)
    \end{equation}
is a Lie subalgebra of $N_{k_1}$, and also 
a Heisenberg 
algebra with a polarization given by the 
first two summands 
in (\ref{polarization}).
We denote the Lie algebra in
(\ref{polarization}) (and also its corresponding Lie
group) by $\overline N_{k_1}$.
The Lie bracket when restricted to the polarization
is described as follows.
\begin{lemma}\label{lemmafortrans}
Let 
$X\in\mathrm{Hom}_D(Y_{k_1+1},X^*_{k_1})$ and
$Y\in\mathrm{Hom}_D(Y_{k_1+1}^*,X^*_{k_1})$. Then 
$[\tilde X,\tilde Y]$ is an element of
$\mathrm{Hom}_D^{\mathrm{inv}}(X_{k_1},X_{k_1}^*)$
given by
$$[\tilde X,\tilde Y]=XY^t-YX^t$$
where $X^t\in\mathrm{Hom}_D(X_{k_1},Y_{k_1+1}^*)$ 
and $Y^t\in\mathrm{Hom}_D(X_{k_1},Y_{k_1+1})$
are uniquely determined as follows
$$\forall\ i\leq k_1\ 
\mathrm{and}\ \forall\ j>k_1,
(X^te_i,e_j)+(e_i,Xe_j)=0$$
$$\forall\ i\leq k_1\ 
\mathrm{and}\ \forall\ j>k_1, 
(Y^te_i,e_j^*)+(e_i,Ye_j^*)=0.$$
\end{lemma}
\begin{proof}
    Follows immediately from 
    $[\tilde X,\tilde Y]=\tilde{X}\tilde{Y}-\tilde{Y}\tilde{X}$ where
    $\tilde{X}$ and $\tilde{Y}$ are defined as in 
    (\ref{transpose}). 
\end{proof}    
    
The adjoint action of 
$\mathrm{Hom}_D^{\mathrm{inv}}
(Y_{k_1+1},Y_{k_1+1}^*)$ on
$\overline N_{k_1}$ is given
by Lemma \ref{actionHOM}.  
This action normalizes $\overline N_{k_1}$ 
and takes 
$\mathrm{Hom}_D(Y_{k_1+1}^*,X_{k_1}^*)$
to
$\mathrm{Hom}_D(Y_{k_1+1},X_{k_1}^*).$\\

At this point we come back to nilpotent 
groups and their representations.
Consider an irreducible 
representation 
$\rho_{1}$ of the Heisenberg
group $N_{k_1}$ with (nontrivial)
central character $\chi_1$.
From the orthogonal decomposition
obtained in (\ref{polarization})
it follows that the restriction of 
$\rho_{1}$ to the group 
$\overline N_{k_1}$
decomposes into 
a direct integral of
representations of this latter Heisenberg
group with the same
central character.
We study the restriction of
a rankable representation of rank one
of the H-tower 
unipotent radical of $G$ to its subgroup
$$\overline N_{k_1}\rtimes 
\mathrm{Hom}_D^{\mathrm{inv}}
(Y_{k_1+1},Y_{k_1+1}^*).$$
This restriction is 
a direct integral of 
representations of the latter group
obtained by extending the irreducible 
representation of $\overline N_{k_1}$
with central character $\chi_1$ to $\overline 
N_{k_1}\rtimes
\mathrm{Hom}_D^{\mathrm{inv}}
(Y_{k_1+1},Y_{k_1+1}^*)$ 
as suggested by Proposition
\ref{weilrepresentationextension}.
This is because
Lemma \ref{commutingspaces} implies that as
subspaces of $N_{k_1}$,
$$[\mathrm{Hom}_D^{\mathrm{inv}}
(Y_{k_1+1},Y_{k_1+1}^*),
\mathrm{Hom}_D(V_c,X^*_{k_1})]=\{0\}
$$
and the Weil representation is functorial. 
(See (1.15) of \cite{Howe} for a precise
meaning of functoriality.)\\

Since 
$\mathrm{Hom}_D^{\mathrm{inv}}(Y_1,Y_1^*)$ 
	is an abelian Lie
group, any unitary
representation of 
this group
can be described as a direct integral of unitary 
characters. 
$\mathrm{Hom}_D^{\mathrm{inv}}(Y_1,Y_1^*)$ 
is
isomorphic to some 
$\mathbb R^p$; so
its group of unitary 
characters
can be naturally 
identified to 
$\mathrm{Hom}_D^{\mathrm{inv}}(Y_1^*,Y_1)$
via the $\mathbb F$-bilinear
form 
\begin{equation}\label{trace}
\beta(A,B)=\mathrm{tr}(AB).
\end{equation}
\begin{definition}
The rank of a character of 
$\mathrm{Hom}_D^{\mathrm{inv}}(Y_1,Y_1^*)$
is  
the $\mathbb F$-rank of the element
in $\mathrm{Hom}_D^{\mathrm{inv}}(Y_1^*,Y_1)$
which corresponds to it via
the bilinear form in (\ref{trace}).
\end{definition}
\begin{definition}\label{Howedefinition} 
    (See \cite{Howe},\cite{LiInv}) 
    Let $\pi$ be a unitary representation
    of $G$. $\pi$ is said to have rank $k$ if and only if
$\pi_{|\mathrm{Hom}_D^{\mathrm{inv}}
	(Y_1,Y_1^*)}$ decomposes into a direct 
	integral of characters of rank equal to $k$.
\end{definition}
The following proposition is a key result
of this section.
\begin{proposition}\label{rankone}
The restriction of a rankable representation
of rank one to 
$$\mathrm{Hom}_D^{\mathrm{inv}}
(Y_1,Y_1^*)$$
is supported on characters whose  
rank is equal to $k_1$.
\end{proposition}

\begin{proof}
The polarization for the group $\overline N_{k_1}$ 
has the structure of a $D$-vector space. 
Therefore, similar to
(\ref{heisenbergaction}), we
can realize the representation
$\rho_{1}$ of
 $\overline N_{k_1}$ on 
 $$L^2(\mathrm{Hom}_D(Y_{k_1+1}^*,X^*_{k_1}))$$
and then extend it to 
$\mathrm{Hom}_D^{\mathrm{inv}}
(Y_{k_1+1},Y_{k_1+1}^*)$.
In Lemma \ref{weilformula} below,
we denote elements of the Lie algebras
by $X,Y,...$ and the corresponding
elements 
in the Lie groups by $e^X,e^Y,$....

\begin{lemma}\label{weilformula}
The action of the extension of $\rho_{1}$
is described as below.
\begin{itemize}
    \item[(a)]For any 
    $X\in\mathrm{Hom}_D(Y_{k_1+1},X_{k_1}^*), 
Y\in\mathrm{Hom}_D(Y^*_{k_1+1},X_{k_1}^*)$
$$(\rho_{1}(e^X)f)(Y)=\chi_1(e^{[Y,X]})f(Y)$$
\item[(b)]
For any $X\in \mathrm{Hom}^{\mathrm{inv}}_D
(Y_{k_1+1},Y_{k_1+1}^*),
Y\in\mathrm{Hom}_D(Y^*_{k_1+1},X_{k_1}^*)$
$$(\rho_{1}(e^X)f)(Y)=\chi_1(e^{{1\over 2}
[Y,YX]})f(Y)$$
\item[(c)] For any $X\in\mathrm{Hom}_D^{\mathrm{inv}}
(X_{k_1},X_{k_1}^*), 
Y\in\mathrm{Hom}_D(Y^*_{k_1+1},X_{k_1}^*)$
$$(\rho_{1}(e^X)f)(Y)=\chi_1(e^X)f(Y).$$
\end{itemize}
\end{lemma}
\begin{proof}
    This is an almost immediate consequence of the 
    Schr\"{o}dinger model for the 
    realization of Weil representation. See 
    \cite{HoweNotes}.
    \end{proof}

Let $X\in\mathrm{Hom}_D^{\mathrm{inv}}
(Y_{k_1+1},Y_{k_1+1}^*)$ and 
$Y\in\mathrm{Hom}_D(Y_{k_1+1}^*,X_{k_1}^*)$. 
For all $i\leq k_1$ and $j>k_1$,
$$
(e_j,(YX)^te_i)=-(YXe_j,e_i)=(Xe_j,Y^te_i)
$$   
$$ =-(e_j,XY^te_i)=(e_j,-XY^te_i)$$   
which implies that $(YX)^t=-XY^t$. 
Thus the equation in
part (b) of
Lemma \ref{weilformula} can be simplified as 
$$(\rho_{1}(X)f)(Y)=\chi_1(e^{-YXY^t})f(Y).$$

We would like to have a single formula instead of 
parts (a), (b) and (c) of Lemma \ref{weilformula}.
To this end, we
define the linear operator $S=S(Y)$ such that
$$S:Y^*_1\mapsto X^*_{k_1}$$ 
by
$$Se^*_i=e^*_i \ \ \   \mathrm{if}\ \ i\leq k_1$$
$$Se^*_i=Ye_i^* \ \ \   \mathrm{if} \ \ i>k_1.$$ 
We have the following lemma.
\begin{lemma}\label{Weilaction}
    Let
$X\in\mathrm{Hom}_D^\mathrm{inv}(Y_1,Y_1^*).$
Let $f$ be a function such that
$$f\in L^2(\mathrm{Hom}_D(Y^*_{k_1+1},X_{k_1}^*)).$$
Then 
$$\rho_{1}(e^X)f(Y)=\chi_1(e^{-SXS^t})f(Y)$$
where $S^t$ is defined as in Lemma \ref{lemmafortrans}.
\end{lemma}
\begin{proof}
    Applying $(e_i,Se_j^*)+(S^te_i,e_j^*)=0$, one can
see that for any $i\leq k_1$
$$S^te_i=-e_i+Y^te_i$$
and thus for any 
$X\in\mathrm{Hom}_D(Y_{k_1+1},X_{k_1}^*)$
\begin{eqnarray*}
-SXS^te_i=-S\tilde{X}(-e_i+Y^te_i)=\\
(-S)(-X^t)e_i-SXY^te_i
=\\
YX^te_i-XY^te_i=[Y,X]e_i
\end{eqnarray*}
which proves $-SXS^t=[Y,X]$.\\

For any $X\in\mathrm{Hom}^{\mathrm{inv}}_D
(Y_{k_1+1},Y_{k_1+1}^*)$ we have
$$-SXS^te_i=-S\tilde{X}(-e_i+Y^te_i)=
-SXY^te_i=-YXY^te_i$$
which proves $-SXS^t=-YXY^t$.\\

Finally, when $X\in\mathrm{Hom}_D^{\mathrm{inv}}
(X_{k_1},X_{k_1}^*)$
$$-SXS^te_i=-SX(-e_i+Y^te_i)=SXe_i=Xe_i$$
which proves $-SXS^t=X$. This completes the proof.
\end{proof}

To complete the proof of Proposition
\ref{rankone}, note that
via the identification described in (\ref{trace}),
the character $\chi(e^X)=\chi_1(e^{-S^tXS})$ 
corresponds to an element 
\begin{equation}\label{svs}
-{a\over k_1}S^tVS\qquad a\in i\mathbb R-\{0\}
\end{equation}
where
$$V:X^*_{k_1}\mapsto X_{k_1}$$
is defined as $$Ve_l^*=(-1)^{l+1}e^*_{k_1-l+1}
\ \ \ \ \mathrm{for\ any}
\ \ \ \ 1\leq l\leq k_1.$$
It is easy to see that (\ref{svs}) is an element of
$\mathrm{Hom}_D^{\mathrm{inv}}(Y_1^*,Y_1)$
of rank $k_1$, even when its domain
is restricted to $X_{k_1}$.
\end{proof}

The following theorem shows that in groups of 
type I, the two notions of rank are essentially the same.
\begin{theorem}\label{maininrank}
    Let $\pi$ be an irreducible
representation of a classical group
    $G$ of type I. Let $N_\gg$ be the H-tower
subgroup of $G$ (see (\ref{levidecompositionofpgamma})).
\begin{itemize}
\item
Assume that the rank of
$\pi$ in the sense of Definition \ref{Howedefinition}
is less than $\hgt(N_\gg)\times k_1$. 
Then 
$\pi$ has rank $k$ in the sense
of Definition \ref{definitionofrank} 
if and only if $\pi$ has rank $kk_1$
in the sense of Definition \ref{Howedefinition}.

\item 
If the rank of $\pi$ in the sense of 
Definition \ref{definitionofrank} is equal
to $\hgt(N_\gg)$, then rank of $\pi$ in the sense
of Definition \ref{Howedefinition} is $\hgt(N_\gg)\times k_1$
or higher.
\end{itemize}

\end{theorem}    

\begin{proof}
    Since $\pi$ is supported on rankable
    representations, 
    the only thing we have to show is that 
   a rankable representation of 
      rank $k$ (in the sense
    of Definition \ref{s-OKPdef}), when restricted to 
    $\mathrm{Hom}_D^{\mathrm{inv}}(Y_1^*,Y_1)$, decomposes
    as a direct integral of 
    characters of rank $kk_1$. By Proposition \ref{rankone},
    this is 
    true when $k=1$.\\
    
    Next consider a rankable 
    representation $\rho$ of rank $k>1$, say 
   $$\rho=\rho_{1}\otimes\cdots\otimes\rho_{k}.$$
    Elementary properties of
    tensor product imply that the restriction of 
    $\rho$ to $\mathrm{Hom}_D^{\mathrm{inv}}(Y_1,Y_1^*)$
    is a direct integral of characters of the form 
    $$\phi_1\cdot\phi_2$$ where the characters
    $\phi_1$ and $\phi_2$ are
    constituents of the direct integral decomposition of 
    $\rho_{1}$ and $\rho_{2}\otimes\cdots\otimes
    \rho_{k}$, when restricted to 
    $\mathrm{Hom}_D^{\mathrm{inv}}(Y_1,Y_1^*)$,
    respectively. 
    But if $\phi_i$ $(i\in\{1,2\})$
    corresponds to 
    $A_i\in
    \mathrm{Hom}_D^{\mathrm{inv}}(Y_1^*,Y_1)$
    via (\ref{trace}), 
    then 
    $\phi_1\cdot\phi_2$ corresponds to
    $A_1+A_2$.
    Since $\rho_{2}\otimes\cdots\otimes
    \rho_{k}$ and therefore any possible
    $\phi_2$ is a trivial representation when
    restricted to $N_{k_1}$, 
    any possible $A_2$ is really an element of 
    $\mathrm{Hom}_D^{\mathrm{inv}}(Y_{k_1+1}^*,Y_{k_1+1})$
    which is extended trivially 
    on $X_{k_1}$ to $Y_1^*$.
   However,  
    at the end of the proof of Proposition \ref{rankone}
    it was shown that any such character $\phi_1$ 
    corresponds to some element
    of $\mathrm{Hom}_D^\mathrm{inv}
    (Y_1^*,Y_1)$ which is of rank $k_1$ even
    when its domain is
    restricted to $X_{k_1}^*$. It is now easy to show that
    we have 
    $$\mathrm{rank}(A_1+A_2)=
    \mathrm{rank}(A_1)+
    \mathrm{rank}(A_2)=k_1+\mathrm{rank}(A_2).$$
    An induction on $k$ completes the proof.
\end{proof}    

\noindent {\bf Example.} Let $G=SO(6,6)$. Let
$\pi$ be an irreducible  unitary representation
of $G$.
Then the rank of $\pi$
in the sense of Definition 
\ref{Howedefinition} can be 0, 2, 4 or 6. The rank of
$\pi$ in the sense of Definition \ref{definitionofrank}
can be 0, 1 or 2. The following chart shows how the
ranks correspond to each other:

\begin{center}
\begin{tabular}{r|l}
Definition \ref{Howedefinition}& Definition \ref{definitionofrank}\\
\hline
0 &0\\
2 & 1\\
4 & 2\\
6&2\\
\end{tabular}
\end{center}

Now let $G=SO(5,11)$. Then we have a similar chart for the rank of $\pi$.

\begin{center}
\begin{tabular}{r|l}
Definition \ref{Howedefinition}& Definition \ref{definitionofrank}\\
\hline
0 &0\\
2 & 1\\
4 & 2\\

\end{tabular}
\end{center}
Therefore the correspondence of ranks may or may
not be one to one. It is an easy exercise to 
determine in which cases the correspondence is
actually one to one.\\

\subsection{$\mathbf{SL_{l+1}}(\mathbb R)$}    
Let $G=SL_{l+1}(\mathbb R)$, the group of
linear transformations on the $l+1$-dimensional
vector space $V$ with a fixed basis 
$$\{e_1,\ldots,e_{l+1}\}.$$
Set
$r=\lfloor {l+1\over 2}\rfloor$.
For any $k$ let $P_k$ be
the maximal parabolic of $G$ which is 
represented
by matrices of the form 
\begin{displaymath}
    \left[
    \begin{array}{cc}
	A&B\\
	0& C\\
    \end{array}
    \right]
\end{displaymath}
where 
$$A\in GL_k(\mathbb R)\ \ ,\ \ 
C\in GL_{l+1-k}(\mathbb R)\ \ , \ \ 
{\rm and}\ \ B\in\mathrm{M}_{k\times (l+1-k)}
(\mathbb R).$$
Therefore 
$$P_k=S(GL_k(\mathbb R)\times GL_{l+1-k}
(\mathbb R))\cdot N_k$$
where $$N_k=\mathrm{Hom}(\mathrm{Span}
(\{e_{k+1},\ldots,
e_{l+1}\}),\mathrm{Span}(\{e_1,\ldots,e_k\})).$$
The parabolic of $G$ which gives the most refined
rank is $P_r$, and 
henceforth we focus our 
attention to $P_r$. Let
\begin{displaymath}
    X_k=\mathrm{Span}(\{e_2,\ldots,e_k\})\qquad,\qquad
Y_k=\mathrm{Span}(\{e_k,\ldots,e_l\}).
\end{displaymath}
The Heisenberg parabolic subgroup of $G$ is $P_1\cap P_l$, 
and a polarization of  
the Lie algebra of its unipotent radical $N$ is 
a direct sum of 
$$\mathrm{Hom}(X_r,\mathbb Re_1)\oplus 
\mathrm{Hom}
(\mathbb Re_{l+1},Y_{r+1})$$
and 
$$\mathrm{Hom}(Y_{r+1},\mathbb Re_1)\oplus
\mathrm{Hom}
(\mathbb Re_{l+1},X_r).
$$
The second summand lies inside the nilradical of
the Lie algebra of $P_r$.
Its center is 
isomorphic to
$$\mathrm{Hom}
(\mathbb R e_{l+1},\mathbb Re_1).$$
As before we are interested in description 
of the restriction of 
a representation $\rho_{1}$ of $N$ to 
$$\mathrm{Hom}
(\mathrm{Span}
(\{e_{r+1},\ldots,
e_{l+1}\}),\mathrm{Span}(\{e_1,\ldots,e_r\})).$$
We identify the dual of 
$$\mathrm{Hom}
(\mathrm{Span}(\{e_{r+1},\ldots,e_{l+1}\}),\mathrm{Span}
    (\{e_1,\ldots,e_r\}))$$
with itself via the bilinear form 
\begin{displaymath}\label{traceform}
    \beta(X,Y)=\mathrm{tr}(X^tY).
    \end{displaymath}
The rank of a unitary character of 
$$\mathrm{Hom}(
\mathrm{Span}
(\{e_{r+1},\ldots,e_{l+1}\}),
\mathrm{Span}
    (\{e_1,\ldots,e_r\}))$$
    is defined to be the rank of the linear transformation
    which corresponds to it via the bilinear form $\beta$.\\
    
    We write any
    $$Y\in\mathrm{Hom}(X_r,\mathbb Re_1)
    \oplus\mathrm{Hom}(\mathbb Re_{l+1},Y_{r+1})$$
    naturally as $Y=Y_1\oplus Y_2$. Define
    $$Y_1^+\in\mathrm{Hom}
    (\mathrm{Span}(\{e_1,\ldots,e_r\}),\mathbb Re_1)$$
    by
    \begin{eqnarray*}
	Y_1^+e_1&=&e_1\\Y_1^+e_j&=&Y_1e_j
	\quad\mathrm{for\ any\ }j>1.
\end{eqnarray*}
Similarly, 
$$Y_2^+\in\mathrm{Hom}(\mathbb Re_{l+1},
\mathrm{Span}(\{e_{r+1},\ldots,e_{l+1}\}))$$
is defined to be 
$$Y_2^+e_{l+1}=-Y_2e_{l+1}+e_{l+1}.$$
We can prove the following version of Lemma
\ref{Weilaction}.
\begin{lemma}
    Let $\rho_{1}$ be a representation 
    of $N$ with central character 
    $\chi_1$ realized on $$\mathcal{H}=L^2
    (\mathrm{Hom}
    (X_r,\mathbb Re_1)\oplus \mathrm{Hom}
(\mathbb Re_{l+1},Y_{r+1}))$$
    as in section 2.1.\\
   For any 
   $$X\in\mathrm{Hom}
    (\mathrm{Span}(\{e_{r+1},\ldots,e_{l+1}\}),\mathrm{Span}
    (\{e_1,\ldots,e_r\})),$$ 
    we have
    $$\rho_{1}(X)f(Y)=\chi_1(e^{Y_1^+XY_2^+})f(Y).$$
    \end{lemma}
\begin{proof}
    We will write $Y$ as  
    $Y_1\oplus Y_2$ according to the polarization given 
    in the statement of the lemma. 
    Based on the Schr\"{o}dinger model, if 
    $X\in \mathrm{Hom}(Y_{r+1},X_r)$, then the 
    action of $X$ on the function $f$ at a point $Y$ is
    multiplication by the character 
    $$\chi_1(e^{-{1\over 2}[Y,YX]})=
    \chi_1(e^{-Y_1XY_2}).$$
    If $X=X_1\oplus X_2\in\mathrm{Hom}(Y_{r+1},
    \mathbb Re_1)\oplus
    \mathrm{Hom}(\mathbb Re_{l+1},X_{r})
    $ 
    then $[X,Y]=-X_1Y_2+Y_1X_2$ and 
    the action of $X$ is by the character
    $$\chi_1(e^{-X_1Y_2+Y_1X_2})$$
and if $X$ belongs to the center of $N$ then
clearly the action will be by the character 
$$\chi_1(e^X).$$
The statement of the lemma follows by a simple
calculation.
\end{proof}    
    
One can see that by the duality provided via 
bilinear form $\beta$, The restriction of a 
representation
of rank one to the nilradical of $P_r$
is a direct integral of characters which correspond
to linear operators of the form
$Y_2^+Y_1^+$, which have
rank one (in the usual sense) 
even when the domain is restricted to 
$\mathbb Re_{l+1}$. Proof of the following theorem 
(which shows the equivalence of the two notions of rank)
is similar to that of Theorem 
\ref{maininrank}.
\begin{theorem}
Let $G=SL_{l+1}(\mathbb R)$.
    Then the restriction of a pure-rank representation of 
    $G$ of rank $k$ (in the sense of 
    Definition \ref{definitionofrank})  
    to the
    abelian nilradical of $P_r$ is supported on 
    unitary characters
    of rank $k$. 
    \end{theorem}    
\noindent {\bf Remark.} Note that when $l$ is even, 
the maximum rank of the unitary characters is $l\over 2$,
which is the same as the height of the H-tower subgoup of $G$.
However, for odd $l$, the maximum rank of unitary 
characters exceeds the height of the H-tower group by one.

\vspace{5mm}


\begin{thebibliography}{ZZZ}
    \bibitem[AV]{adamsvogan}
    Adams, J. and Vogan, D. (Editors)
    {\it Representation theory of Lie groups},
    Lectures from the Graduate Summer School
    held in Park City, UT, July 13--31, 1998. 
   IAS/Park City Mathematics Series, 8, 2000. 
   \bibitem[Ad]{AdamsTrapa}
    Adams, J. 
    {\it Non-linear covers of real groups}, 
preprint, 2004.
\bibitem[Bor]{Borel} Borel, A., 
{\it Linear algebraic groups}, Springer-Verlag, 
New York, 1991.
\bibitem[Bou]{bourbaki} Bourbaki, N., 
{\it Lie groups and Lie algebras. Chapters 4--6.}
 Translated from the 1968 French original by Andrew
Pressley. Elements of Mathematics (Berlin). 
Springer-Verlag, Berlin, 2002.





\bibitem[Br]{brown} Brown, I. D., {\it Dual topology of a nilpotent 
Lie group}, Ann. Sci. de l'E.N.S. $4^e$ serie, no. 6, (1973), pp 407-411.
   \bibitem[CG]{CG}
    Corwin, L. J. and Greenleaf, F. P., 
    {\it Representations of nilpotent Lie groups and 
    their applications}, Cambridge university press, 
    1990.
    \bibitem[De]{deod}Deodhar, Vinay V.
    {\it On central extensions of rational points of algebraic groups},
    Amer. J. Math. 100 (1978), no. 2, 303--386.
\bibitem[Du]{dufloacta} Duflo, M., {\it Th\'{e}orie de Mackey pour les groupes de Lie 
alg\'{e}brique}, Acta Math. 149 (1982), 152--213.



\bibitem[DS]{SAHDEV}
Dvorsky, A and Sahi, S., {\it Tensor products of singular 
representations and an extension of the $\theta$-correspondence}. 
Selecta Math. (N.S.) 4 (1998),
no. 1, 11--29.









    \bibitem[GS]{GanSavin}
Gan, W. T. and Savin, Gordan {\it On minimal representations},
Representation Theory 9 (2005), 46--93.
    
\bibitem[GRS]{Ginz}
Ginzburg, David and  Rallis, Stephen and Soudry, David,
{\it A tower of theta correspondences for $G\sb 2$}. 
Duke Math. J. 88 (1997), no. 3, 537--624.


















    \bibitem[GW]{GW} Gross, B. H. and Wallach, N. R.,
    {\it On quaternionic discrete series
    representations, and their continuations},
    J. Reine Angew. Math. 481 (1996), 73--123.
    \bibitem[Hel]{Helgason} Helgason, S.
    {\it Differential geometry, Lie groups, and
    symmetric spaces}, AMS, 2000.
    \bibitem[Hw1]{Howe} Howe, R.,
    {\it On a notion of rank for unitary
    representations of the classical groups},
    Harmonic analysis and group representations, 223--331,
    Liguori, Naples, 1982.
    \bibitem[Hw2]{HoweBulletin}
   Howe, R.
   {\it On the role of the Heisenberg group
   in harmonic analysis.}, Bull. Amer. Math. Soc.
   (N.S.) 3 (1980), no. 2, 821--843.
    \bibitem[Hw3]{howeconstruct} Howe, R.
    {\it Small unitary representations of classical groups},
    Group representations, ergodic theory,
    operator algebras, and mathematical physics
    (Berkeley, Calif., 1984), 121--150,
    Math. Sci. Res. Inst. Publ., 6,
    Springer, New York, 1987.
    \bibitem[Hw4]{Howeautomorphic} Howe, R.,
    {\it Automorphic forms of low rank},
    Noncommutative harmonic analysis
    and Lie groups (Marseille, 1980),
    pp. 211--248, Lecture Notes in Math.,
    880, Springer, Berlin-New York, 1981.
    \bibitem[Hw5]{Howetheta}
    Howe, R., {\it $\theta$-series and 
invariant theory},
    Automorphic forms,
    representations and $L$-functions
    (Proc. Sympos. Pure Math., Oregon
    State Univ., Corvallis, Ore., 1977),
    Part 1, pp. 275--285
    \bibitem[Hw6]{HoweNotes} Howe, R., {\it Oscillator 
    representation, analytic preliminaries}, 
    unpublished notes, 1977.
    \bibitem[HM]{HoweMoore} Howe, R. and Moore, C.
    {\it Asymptotic properties of unitary representations}, 
    J. Funct. Anal. 32 (1979), no. 1, 72--96.
    \bibitem[HRW]{HoweOKP}Howe, R., Ratcliff, G. and 
    Wildberger, N.,
    {\it Symbol mappings for certain nilpotent groups}, 
    Lie group representations, III (College Park, Md., 1982/1983), 288--320, 
    Lecture Notes in Math., 1077, 
    Springer, Berlin, 1984. 
    
\bibitem[Hu]{humph} Humphreys, J. E., 
{\it Linear algebraic groups}, Graduate 
Texts in Mathematics 21, Springer-Verlag,
New York-Heidelberg, 1975.




\bibitem[Jo]{joseph}
    Joseph, A. {\it 
    A preparation theorem for the prime 
    spectrum of a semisimple Lie algebra},
    J. Algebra 48 (1977), no. 2, 241--289.
\bibitem[Kaz]{kazhdan}
Kazhdan, David, {\it 
Some applications of the Weil representation.} 
J. Analyse Mat. 32 (1977), 235--248.



     \bibitem[Kr1]{Kirillov}
    Kirillov, A. A., 
    {\it Unitary representations of nilpotent Lie groups}
    (Russian), 
    Uspehi Mat. Nauk 17 1962 no. 4 (106), 57--110.
    \bibitem[Kr2]{Kirillov2}
Kirillov, A. A., {\it Lectures on the orbit method},
Graduate Studies in Mathematics 64, American 
Mathematical Society, Providence, RI 2004.




\bibitem[Ki]{AlexanderKirillov} Kirillov, A. A. Jr.,
    {\it Merits and demerits of the orbit method}, 
    (English. English summary) 
    Bull. Amer. Math. Soc. (N.S.) 36 (1999), no. 4, 433--488.
    \bibitem[Li1]{LiInv} Li, J.-S.,
    {\it Singular unitary representations of 
    classical groups},
    Invent. Math. 97 (1989), no. 2, 237--255.
    \bibitem[Li2]{LiComp}	
    Li, J. -S.,
    {\it On the classification of irreducible 
    low rank unitary representations of classical groups},
    Compositio Math. 71 (1989), no. 1, 29--48.
    \bibitem[Li3]{Liautomorphic}
    Li, J. -S., {\it Automorphic forms with 
    degenerate Fourier coefficients},
    Amer. J. Math. 119 (1997), no. 3, 523--578.
    
\bibitem[LZ]{LiZhu}
Li, J. S., Zhu, C. B.,
{\it On the decay of matrix coefficients 
for exceptional groups},
Math. Ann. {\bf 305},  no. 2, 249--270 (1996).








\bibitem[LS]{savinloke} Loke, H. Y. and Savin, G., 
{\it Rank and matrix
    coefficents for simply-laced groups}, preprint, 2005.
    \bibitem[Mac1]{Mackey}
    Mackey, G. W., {\it The theory of unitary group
    representations}, University of Chicago Press, 1976.
\bibitem[Mac2]{MackeyI} Mackey, G. W.,
{\it Induced representations of locally
compact groups I},  Ann. of Math. (2) 55, 
(1952). 101--139.

\bibitem[MS]{MAG}
Magaard, K. and Savin, G., {\it  
Exceptional $\Theta$-correspondences. I}. Compositio Math. 107 (1997), 
no. 1, 89--123.










\bibitem[Mo1]{moore}
    Moore, Calvin C., {\it Decomposition of unitary representations
    defined by discrete subgroups of nilpotent groups},
    Ann. of Math. (2)  82  1965 ,146--182.

\bibitem[Mo2]{mooreextension}
Moore, Calvin C., 
{\it Extensions and low dimensional cohomology 
theory of locally compact groups I, II.}
Trans. Amer. Math.
  Soc. 113 1964 40--63.








\bibitem[Oh]{Oh}

Oh, H., {\it Uniform pointwise bounds 
for matrix coefficients of 
unitary representations and applications to 
Kazhdan constants},
Duke Math. J.  {\bf 113 } no. 1, 133--192 (2002). 









\bibitem[PR]{plra}
Platonov, V. and Rapinchuk, A.,
{\it 
Algebraic groups and number theory}. 
Translated from the 1991 Russian original 
by Rachel Rowen. Pure and Applied Mathematics, 139. 
Academic Press, Inc., Boston, MA, 1994.    







\bibitem[Sa]{hadi2}
    Salmasian, H., {\it Isolatedness of minimal representations and
    minimal decay of matrix coefficients}, submitted.
    \bibitem[Sca]{Scara}
    Scaramuzzi, R.,
    {\it A notion of rank for
    unitary representations of general linear groups},
    Trans. Amer. Math. Soc. 319 (1990),
    no. 1, 349--379.
	\bibitem[Se]{serre} Serre, J. P., {\it
Lie algebras and Lie groups}, Lecture Notes 
in Mathematics 1500, Springer, Berlin, 1992.


    \bibitem[Tay]{taylor}
    Taylor, M. E., {\it 
Noncommutative harmonic analysis},
    AMS Mathematical Surveys and 
Monographs, no. 22, 1986.
    \bibitem[Tor]{torassoduke}Torasso, P.
    {\it M\'{e}thode des orbites de Kirillov-Duflo et
    repr\'{e}sentations minimales des groupes
    simples sur un corps local de caract\'{e}ristique nulle}
    (French),
    Duke Math. J. 90 (1997), no. 2, 261--377.
    \bibitem[Ti]{Tits} Tits, Jacques
    {\it Tabellen zu den einfachen Lie
    Gruppen und ihren Darstellungen} (German),
Lecture Notes in Mathematics 40,
    Springer-Verlag, Berlin-New York 1967.
\bibitem[We1]{WEILNUM}
Weil, A., {\it  Basic number theory.} Third edition. 
Die Grundlehren der Mathematischen Wissenschaften, Band 144. 
Springer-Verlag, New York-Berlin, 1974. 

    \bibitem[We2]{weil} Weil, A. {\it Sur certains groupes
    d'op\'{e}rateurs unitaires}, Acta Math. 111 (1964),
    143-211.
    \bibitem[Ws]{Marty} Weissman, M.H., 
    {\it The Fourier-Jacobi map and small 
    representations}, 
Represent. Theory 7 (2003), 275--299.
\end{thebibliography}
\end{document}